\documentclass[reqno,11pt]{amsart}
\usepackage{amsmath}
\usepackage{amsxtra}
\usepackage{amscd}
\usepackage{amsthm}
\usepackage{amsfonts}
\usepackage{amssymb}
\usepackage{eucal}
\usepackage[matrix,arrow,curve]{xy}
\theoremstyle{plain}
\newtheorem{Thm}[subsection]{Theorem}
\newtheorem{Cor}[subsection]{Corollary}
\newtheorem{Lem}[subsection]{Lemma}
\newtheorem{Prop}[subsection]{Proposition}
\newtheorem{Conj}[subsection]{Conjecture}

\theoremstyle{definition}
\newtheorem{Def}[subsection]{Definition}

\theoremstyle{remark}

\newtheorem{Rem}[subsection]{Remark}

\newtheorem{conj}{Conjecture}

\numberwithin{equation}{section}
\renewcommand{\rm}{\normalshape}

\newif\ifShowLabels
\ShowLabelstrue
\newdimen\theight
\def\TeXref#1{%
    \leavevmode\vadjust{\setbox0=\hbox{{\tt
        \quad\quad  {\small \rm #1}}}%
    \theight=\ht0
    \advance\theight by \lineskip
    \kern -\theight \vbox to
    \theight{\rightline{\rlap{\box0}}%
    \vss}%
    }}%

\ShowLabelsfalse% comment this out if labels should be printed

\renewcommand{\sec}[2]{\section{#2}\label{S:#1}%
    \ifShowLabels \TeXref{{S:#1}} \fi}
\newcommand{\ssec}[2]{\subsection{#2}\label{SS:#1}%
    \ifShowLabels \TeXref{{SS:#1}} \fi}

\newcommand{\refs}[1]{Section ~\ref{S:#1}}
\newcommand{\refss}[1]{Section ~\ref{SS:#1}}

\newcommand{\reft}[1]{Theorem ~\ref{T:#1}}

\newcommand{\refp}[1]{Proposition ~\ref{P:#1}}
\newcommand{\refc}[1]{Corollary ~\ref{C:#1}}

\newcommand{\refe}[1]{\eqref{E:#1}}
\newcommand{\refco}[1]{Conjecture ~\ref{Co:#1}}

\newenvironment{thm}[1]%
    { \begin{Thm} \label{T:#1}  \ifShowLabels \TeXref{T:#1} \fi }%
    { \end{Thm} }

\renewcommand{\th}[1]{\begin{thm}{#1} \sl }
\renewcommand{\eth}{\end{thm} }

\newenvironment{lemma}[1]%
    { \begin{Lem} \label{L:#1}  \ifShowLabels \TeXref{L:#1} \fi }%
    { \end{Lem} }
\newcommand{\lem}[1]{\begin{lemma}{#1} \sl}
\newcommand{\elem}{\end{lemma}}

\newenvironment{propos}[1]%
    { \begin{Prop} \label{P:#1}  \ifShowLabels \TeXref{P:#1} \fi }%
    { \end{Prop} }
\newcommand{\prop}[1]{\begin{propos}{#1}\sl }
\newcommand{\eprop}{\end{propos}}

\newenvironment{corol}[1]%
    { \begin{Cor} \label{C:#1}  \ifShowLabels \TeXref{C:#1} \fi }%
    { \end{Cor} }
\newcommand{\cor}[1]{\begin{corol}{#1} \sl }
\newcommand{\ecor}{\end{corol}}

\newenvironment{defeni}[1]%
    { \begin{Def} \label{D:#1}  \ifShowLabels \TeXref{D:#1} \fi }%
    { \end{Def} }
\newcommand{\defe}[1]{\begin{defeni}{#1} \sl }
\newcommand{\edefe}{\end{defeni}}

\newenvironment{remark}[1]%
    { \begin{Rem} \label{R:#1}  \ifShowLabels \TeXref{R:#1} \fi }%
    { \end{Rem} }
\newcommand{\rem}[1]{\begin{remark}{#1}}
\newcommand{\erem}{\end{remark}}

\newenvironment{conjec}[1]%
    { \begin{Conj} \label{Co:#1}  \ifShowLabels \TeXref{Co:#1} \fi }%
    { \end{Conj} }
\renewcommand{\conj}[1]{\begin{conjec}{#1} \sl }
\newcommand{\econj}{\end{conjec}}

\newcommand{\eq}[1]%
    { \ifShowLabels \TeXref{E:#1} \fi
       \begin{equation} \label{E:#1} }
\newcommand{\eeq}{ \end{equation} }

\newcommand{\prf}{ \begin{proof} }
\newcommand{\epr}{ \end{proof} }

\newcommand\alp{\alpha}     

\newcommand\gam{\gamma}     
     \newcommand\Del{\Delta}
\newcommand\eps{\varepsilon}

\newcommand\kap{\kappa}
\newcommand\lam{\lambda}        \newcommand\Lam{\Lambda}

\newcommand\calA{{\mathcal{A}}}

\newcommand\calF{{\mathcal{F}}}
\newcommand\calG{{\mathcal{G}}}

\newcommand\calK{{\mathcal{K}}}

\newcommand\calM{{\mathcal{M}}}

\newcommand\calU{{\mathcal{U}}}

\newcommand\calW{{\mathcal{W}}}

\newcommand\calZ{{\mathcal{Z}}}

\newcommand\bfc{{\mathbf c}}        \newcommand\bfC{{\mathbf C}}

\newcommand\bfg{{\mathbf g}}

        \newcommand\bfM{{\mathbf M}}

        \newcommand\bfX{{\mathbf X}}

\newcommand\TT{\mathbb{T}}

\newcommand\PP{\mathbb{P}}
\renewcommand\AA{\mathbb{A}}

\newcommand\GG{\mathbb{G}}

\newcommand\ZZ{\mathbb{Z}}

\newcommand\CC{\mathbb{C}}

 \newcommand\grg{{\mathfrak{g}}}
 \newcommand\grh{{\mathfrak{h}}}

 \newcommand\grl{{\mathfrak{l}}}
 \newcommand\grm{{\mathfrak{m}}}
 \newcommand\grn{{\mathfrak{n}}}
 
 \newcommand\grp{{\mathfrak{p}}}
 \newcommand\grq{{\mathfrak{q}}}

 \newcommand\grt{{\mathfrak{t}}}

\newcommand\sdp{\times \hskip -0.3em {\raise 0.3ex
\hbox{$\scriptscriptstyle |$}}} % semidirect product

\newcommand\End{\operatorname{End\,}}

\newcommand\GL{\operatorname{GL}}
\newcommand\Gr{\operatorname{Gr}}

\newcommand\Hom{\operatorname {Hom}}

\newcommand\Int{\operatorname{Int}}

\newcommand\SL{{\rm SL}}

\newcommand\Sym{\operatorname{Sym}}

\newcommand\x{\times}
\newcommand\ten{\otimes}

\newcommand{\ra}{\rangle}
\newcommand{\la}{\langle}

\newcommand\nc{\newcommand}

\newcommand{\iso}{{\stackrel{\sim}{\longrightarrow}}}

\nc\aff{\operatorname{aff}}
\nc\oGr{\overline{\Gr}}
\nc\Bun{\operatorname{Bun}}
\nc\hgrg{\widehat{\grg}}
\renewcommand\Int{\operatorname{Int}}
\nc\bInt{\overline{\Int}}
\nc\hatLam{\widehat{\Lam}}
\nc\bmu{\overline{\mu}}
\nc\bnu{\overline{\nu}}
\nc\blambda{\overline{\lam}}
\renewcommand\SL{\operatorname{SL}}
\nc\ocalW{\overline{\calW}}
\nc\pos{\operatorname{pos}}
\nc\IH{\operatorname{IH}}
\nc\Rep{\operatorname{Rep}}
\nc\Gal{\operatorname{Gal}}
\nc{\tilGr}{\widetilde{\Gr}}

\nc\Pic{\operatorname{Pic}}
\nc{\HC}{{\mathcal{HC}}}
\nc{\on}{\operatorname}
\nc{\BA}{{\mathbb{A}}}
\nc{\BC}{{\mathbb{C}}}
\nc{\BM}{{\mathbb{M}}}
\nc{\BN}{{\mathbb{N}}}
\nc{\BP}{{\mathbb{P}}}
\nc{\BR}{{\mathbb{R}}}
\nc{\BZ}{{\mathbb{Z}}}
\nc{\BS}{{\mathbb{S}}}

\nc{\CA}{{\mathcal{A}}}
\nc{\CB}{{\mathcal{B}}}
\nc{\CalD}{{\mathcal D}}
\nc{\CE}{{\mathcal{E}}}
\nc{\CF}{{\mathcal{F}}}
\nc{\CG}{{\mathcal{G}}}
\nc{\CK}{{\mathcal{K}}}
\nc{\CL}{{\mathcal{L}}}
\nc{\CM}{{\mathcal{M}}}
\nc{\CMM}{{\mathcal{M}^{\operatorname{gen}}_\hbar(-\rho)}}
\nc{\CN}{{\mathcal{N}}}
\nc{\CO}{{\mathcal{O}}}
\nc{\CP}{{\mathcal{P}}}
\nc{\CQ}{{\mathcal{Q}}}
\nc{\CR}{{\mathcal{R}}}
\nc{\CS}{{\mathcal{S}}}
\nc{\CT}{{\mathcal{T}}}
\nc{\CU}{{\mathcal{U}}}
\nc{\CV}{{\mathcal{V}}}
\nc{\CW}{{\mathcal{W}}}
\nc{\CX}{{\mathcal{X}}}
\nc{\CZ}{{\mathcal{Z}}}

\nc{\gen}{{\operatorname{gen}}}
\nc{\cM}{{\check{\mathcal M}}{}}
\nc{\csM}{{\check{\mathcal A}}{}}
\nc{\obM}{{\overset{\circ}{\mathbf M}}{}}
\nc{\oCA}{{\overset{\circ}{\mathcal A}}{}}
\nc{\obA}{{\overset{\circ}{\mathbf A}}{}}
\nc{\ooM}{{\overset{\circ}{M}}{}}
\nc{\osM}{{\overset{\circ}{\mathsf M}}{}}
\nc{\vM}{{\overset{\bullet}{\mathcal M}}{}}
\nc{\nM}{{\underset{\bullet}{\mathcal M}}{}}
\nc{\obD}{{\overset{\circ}{\mathbf D}}{}}
\nc{\cp}{{\overset{\circ}{\mathbf p}}{}}
\nc{\ofZ}{{\overset{\circ}{\mathfrak Z}}{}}

\nc{\fa}{{\mathfrak{a}}}
\nc{\fb}{{\mathfrak{b}}}
\nc{\fe}{{\mathfrak{e}}}
\nc{\ff}{{\mathfrak{f}}}
\nc{\fg}{{\mathfrak{g}}}
\nc{\fgl}{{\mathfrak{gl}}}
\nc{\fh}{{\mathfrak{h}}}
\nc{\fj}{{\mathfrak{j}}}
\nc{\fl}{{\mathfrak{l}}}
\nc{\fm}{{\mathfrak{m}}}
\nc{\fn}{{\mathfrak{n}}}
\nc{\fu}{{\mathfrak{u}}}
\nc{\fp}{{\mathfrak{p}}}
\nc{\fq}{{\mathfrak{q}}}
\nc{\frr}{{\mathfrak{r}}}
\nc{\fs}{{\mathfrak{s}}}
\nc{\ft}{{\mathfrak{t}}}
\nc{\fw}{{\mathfrak{w}}}
\nc{\fT}{{\mathfrak{T}}}
\nc{\ofT}{{\overline{\mathfrak T}}}
\nc{\ofS}{{\overline{\mathfrak S}}}
\nc{\fsl}{{\mathfrak{sl}}}
\nc{\hsl}{{\widehat{\mathfrak{sl}}}}
\nc{\hgl}{{\widehat{\mathfrak{gl}}}}
\nc{\hg}{{\widehat{\mathfrak{g}}}}
\nc{\chg}{{\widehat{\mathfrak{g}}}{}^\vee}
\nc{\hn}{{\widehat{\mathfrak{n}}}}
\nc{\chn}{{\widehat{\mathfrak{n}}}{}^\vee}

\nc{\fA}{{\mathfrak{A}}}
\nc{\fB}{{\mathfrak{B}}}
\nc{\fD}{{\mathfrak{D}}}
\nc{\fE}{{\mathfrak{E}}}
\nc{\fF}{{\mathfrak{F}}}
\nc{\fG}{{\mathfrak{G}}}
\nc{\fK}{{\mathfrak{K}}}
\nc{\fL}{{\mathfrak{L}}}
\nc{\fM}{{\mathfrak{M}}}
\nc{\fN}{{\mathfrak{N}}}
\nc{\frP}{{\mathfrak{P}}}
\nc{\fQ}{{\mathfrak{Q}}}
\nc{\fS}{{\mathfrak S}}
\nc{\fU}{{\mathfrak{U}}}
\nc{\fV}{{\mathfrak{V}}}
\nc{\fZ}{{\mathfrak{Z}}}

\nc{\ba}{{\mathbf{a}}}
\nc{\bb}{{\mathbf{b}}}
\nc{\bc}{{\mathbf{c}}}
\nc{\be}{{\mathbf{e}}}
\nc{\bj}{{\mathbf{j}}}
\nc{\bn}{{\mathbf{n}}}
\nc{\bp}{{\mathbf{p}}}
\nc{\bq}{{\mathbf{q}}}
\nc{\br}{{\mathbf{r}}}
\nc{\bv}{{\mathbf{v}}}
\nc{\bx}{{\mathbf{x}}}
\nc{\by}{{\mathbf{y}}}
\nc{\bw}{{\mathbf{w}}}
\nc{\bA}{{\mathbf{A}}}
\nc{\bB}{{\mathbf{B}}}
\nc{\bC}{{\mathbf{C}}}
\nc{\bK}{{\mathbf{K}}}
\nc{\bD}{{\mathbf{D}}}
\nc{\bH}{{\mathbf{H}}}
\nc{\bM}{{\mathbf{M}}}
\nc{\bN}{{\mathbf{N}}}
\nc{\bS}{{\mathbf{S}}}
\nc{\bT}{{\mathbf{T}}}
\nc{\bV}{{\mathbf{V}}}
\nc{\bW}{{\mathbf{W}}}
\nc{\bX}{{\mathbf{X}}}
\nc{\bP}{{\mathbf{P}}}
\nc{\bZ}{{\mathbf{Z}}}

\nc{\sA}{{\mathsf{A}}} \nc{\sB}{{\mathsf{B}}} \nc{\sC}{{\mathsf{C}}}
\nc{\sD}{{\mathsf{D}}} \nc{\sF}{{\mathsf{F}}} \nc{\sK}{{\mathsf{K}}}
\nc{\sM}{{\mathsf{M}}} \nc{\sO}{{\mathsf{O}}} \nc{\sQ}{{\mathsf{Q}}}
\nc{\sP}{{\mathsf{P}}} \nc{\sT}{{\mathsf{T}}} \nc{\sV}{{\mathsf{V}}}
\nc{\sZ}{{\mathsf{Z}}} \nc{\sa}{{\mathsf{a}}} \nc{\se}{{\mathsf{e}}}
\nc{\ssf}{{\mathsf{f}}} \nc{\sfp}{{\mathsf{p}}}
\nc{\sr}{{\mathsf{r}}} \nc{\sfb}{{\mathsf{b}}}
\nc{\sfc}{{\mathsf{c}}} \nc{\sd}{{\mathsf{d}}}
\nc{\sfl}{{\mathsf{l}}}

\nc{\BK}{{\bar{K}}}

\nc{\tA}{{\widetilde{\mathbf{A}}}}
\nc{\tB}{{\widetilde{\mathcal{B}}}}
\nc{\tg}{{\widetilde{\mathfrak{g}}}}
\nc{\tG}{{\widetilde{G}}}
\nc{\TM}{{\widetilde{\mathbb{M}}}{}}
\nc{\tO}{{\widetilde{\mathsf{O}}}{}}
\nc{\tU}{{\widetilde{\mathfrak{U}}}{}}
\nc{\TZ}{{\tilde{Z}}}
\nc{\tx}{{\tilde{x}}}
\nc{\tbv}{{\tilde{\bv}}}
\nc{\tfP}{{\widetilde{\mathfrak{P}}}{}}
\nc{\tz}{{\tilde{\zeta}}}
\nc{\tmu}{{\tilde{\mu}}}

\nc{\td}{\widetilde{\underline{d}}{}}
\nc{\tzeta}{\widetilde{\zeta}{}}
\nc{\hd}{{\widehat{\underline{d}}}}
\nc{\hG}{{\widehat{G}}}
\nc{\hBP}{\widehat{\mathbb P}{}}
\nc{\hQ}{{\widehat{Q}}}
\nc{\hsM}{\widehat{\mathsf M}{}}
\nc{\hfM}{\widehat{\mathfrak M}{}}
\nc{\hCP}{\widehat{\mathcal P}{}}
\nc{\hCR}{\widehat{\mathcal R}{}}
\nc{\hCS}{{\widehat{\mathcal S}}}
\nc{\hfZ}{\widehat{\mathfrak Z}{}}

\nc{\urho}{\underline{\rho}}
\nc{\uB}{\underline{B}}
\nc{\uC}{{\underline{\mathbb{C}}}}
\nc{\ui}{\underline{i}}
\nc{\ofP}{{\overline{\mathfrak{P}}}}
\nc{\hrho}{{\hat{\rho}}}

\nc{\unl}{\underline}
\nc{\ol}{\overline}
\nc{\one}{{\mathbf{1}}}
\nc{\two}{{\mathbf{t}}}

\nc{\Tot}{{\mathop{\operatorname{\rm Tot}}}}
\nc{\Hilb}{{\mathop{\operatorname{\rm Hilb}}}}
\nc{\CHom}{{\mathop{\operatorname{{\mathcal{H}}\it om}}}}
\nc{\defi}{{\mathop{\operatorname{\rm def}}}}
\nc{\length}{{\mathop{\operatorname{\rm length}}}}
\nc{\Cliff}{{\mathsf{Cliff}}}
\nc{\Fl}{{\mathsf{Fl}}}
\nc{\Fib}{{\mathsf{Fib}}}
\nc{\Coh}{{\mathsf{Coh}}}
\nc{\FCoh}{{\mathsf{FCoh}}}

\nc{\reg}{{\text{\rm reg}}}

\nc{\cplus}{{\mathbf{C}_+}}
\nc{\cminus}{{\mathbf{C}_-}}
\nc{\cthree}{{\mathbf{C}_*}}
\nc{\Qbar}{{\bar{Q}}}

\nc{\bh}{{\bar{h}}}
\nc{\bOmega}{{\overline{\Omega}}}
\nc\tGr{\widetilde{\Gr}}

\nc{\seq}[1]{\stackrel{#1}{\sim}}

\nc\uS{\underline{S}}
\newcommand\F{\calF}

\newcommand\QM{\mathcal QM}
\nc{\chH}{\check H}
\nc{\chh}{\check \grh}

\renewcommand\chn{\check \grn}

\renewcommand\chg{\check \grg}

\nc\chT{\check T}

\nc\Vir{\mathbf{Vir}}
\nc\chL{\check L}

\begin{document}
\title{A finite analog of the AGT relation I: finite $W$-algebras and quasimaps' spaces}
\author{Alexander Braverman, Boris Feigin, Michael Finkelberg and
Leonid Rybnikov}

\begin{abstract}Recently Alday, Gaiotto and Tachikawa proposed a conjecture
relating 4-dimensional
super-symmetric gauge theory for a gauge group $G$ with certain 2-dimensional
conformal field theory. This conjecture
implies the existence of certain structures on the (equivariant) intersection cohomology of the
Uhlenbeck partial compactification of the moduli space of framed $G$-bundles on $\PP^2$.
More precisely, it predicts the existence of an action of the corresponding $W$-algebra on the above cohomology,
satisfying certain properties.

We propose a ``finite analog" of the (above corollary of the) AGT conjecture.
Namely, we replace
the Uhlenbeck space with the space of {\em based quasi-maps} from $\PP^1$ to any partial flag variety
$G/P$ of $G$ and conjecture that its equivariant intersection cohomology carries an action of the finite
$W$-algebra $U(\grg,e)$ associated with the principal nilpotent element in the
Lie algebra of the Levi subgroup of $P$; this action
is expected to satisfy some list of natural properties. This conjecture
generalizes the main result of \cite{B} when $P$ is the Borel subgroup. We prove our conjecture for
$G=GL(N)$, using the works of Brundan and Kleshchev interpreting the algebra $U(\grg,e)$ in terms of certain
shifted Yangians.
\end{abstract}
\maketitle
%%%%%%%%%%%%%%%%%%%%%%%%%%%%%%%%%%%%%%%%%%%%%%%%%%%%%%%%%%%%%%%%%%%%%%%%%%%%%%%%%%%%%%%%%%%%%%%%%%%%%%%%%%%%%%%%
\sec{int}{Introduction}
%%%%%%%%%%%%%%%%%%%%%%%%%%%%%%%%%%%%%%%%%%%%%%%%%%%%%%%%%%%%%%%%%%%%%%%%%%%%%%%%%%%%%%%%%%%%%%%%%%%%%%%%%%%%%%%%%%%

\ssec{setup}{The setup}Let $G$ be a semi-simple simply connected complex algebraic group (or, more generally, a connected reductive
group whose derived group $[G,G]$ is simply connected) and let
$P$ be a parabolic subgroup of $G$. We shall denote by $L$ the corresponding Levi
factor. Let $B$ be a Borel subgroup of $G$ contained in $P$ and containing a maximal
torus $T$ of $G$. We shall also denote by $\Lam$ the coweight lattice of
$G$ (which is the same as the lattice of cocharacters of $T$); it has a quotient lattice
$\Lam_{G,P}=\Hom(\CC^*,L/[L,L])$, which can also be regarded as the lattice of characters
of the center $Z(\check L)$ of the Langlands dual group $\check L$.
Note that $\Lam_{G,B}=\Lam$. The lattice $\Lam_{G,P}$ contains
canonical sub-semi-group $\Lam_{G,P}^+$ spanned by the images of positive coroots of $G$.

Set now $\calG_{G,P}=G/P$. There is a natural isomorphism $H_2(\calG_{G,P},\ZZ)\simeq \Lam_{G,P}$.
Let $C$ be a smooth connected projective curve over $\CC$.
Then the degree of a map $f:C\to \calG_{G,P}$ can
be considered as an element of $\Lam_{G,P}$; it is easy to see that actually $\deg f$ must lie
in $\Lam_{G,P}^+$. For $\theta,\theta'\in\Lam_{G,P}$ we shall write $\theta\geq \theta'$ if
$\theta-\theta'\in\Lam_{G,P}^+$.

Let
$e_{{G,P}}\in\calG_{G,P}$ denote the image of $e\in G$.
Clearly $e_{G,P}$ is stable under the action of $P$ on
$\calG_{G,P}$. Let $\calM_{G,P}$ denote the moduli space of based
maps from $(\PP^1,\infty)$ to
$(\calG_{G,P},e_{{G,P}})$, i.e. the moduli space of maps
$\PP^1\to\calG_{G,P}$ which send $\infty$ to $e_{G,P}$.
This space is acted on by the group $P\x\CC^*$, where $P$ acts on $\calG_{G,P}$ preserving the point
$e_{G,P}$ and $\CC^*$ acts on $\PP^1$ preserving $\infty$; in particular, the reductive group
$L\x\CC^*$ acts on $\calM_{G,P}$.
Also for any $\theta\in\Lam_{G,P}^+$ let $\calM_{G,P}^{\theta}$ denote the space of maps as above
of degree $\theta$.

For each
$\theta$ as above one can also consider the space of {\it based
quasi-maps} (or {\it Zastava space} in the terminology of
\cite{BFGM}, \cite{FM} and \cite{FFKM}; cf. also \cite{Bicm} for a review of quasi-maps' spaces) which we denote by
$\QM^{\theta}_{G,P}$. This is an affine algebraic variety
containing $\calM^{\theta}_{G,P}$ as a dense open subset. Moreover, it possesses a stratification of the
form
$$
\QM_{G,P}^{\theta}=\bigcup\limits_{0\leq\theta'\leq\theta}\calM_{G,P}^{\theta'}\x\Sym^{\theta-\theta'}\AA^1,
$$
where for any $\gam\in \Lam_{G,P}^+$ we denote by $\Sym^{\gam}\AA^1$ the variety of formal linear
combinations $\sum \lam_i x_i$ where $x_i\in \AA^1$ and $\lam_i\in \Lam_{G,P}^+$ such that
$\sum\lam_i=\gam$. In most cases the variety $\QM_{G,P}^{\theta}$ is singular.
%-----------------------------------------------------------------------------------------------------
\ssec{}{Equivariant integration}
For a connected reductive group $\GG$ with a maximal torus $\TT$ let
$$
\calA_{\GG}=H^*_{\GG}(pt,\CC).
$$
This is a graded algebra which is known to be canonically isomorphic to the algebra of $\GG$-invariant
polynomial functions on the Lie algebra $\bfg$ of $\GG$. We shall
denote by $\calK_{\GG}$ its field of fractions.
Let now $Y$ be a variety endowed with an action of $\GG$ such that
$Y^{\TT}$ is proper. Then as was remarked e.g. in \cite{B} we have a well-defined
integration map
$$
\int_Y:\IH^*_{\GG}(Y)\to\calK_{\GG},
$$
which is a map of $\calA_{\GG}$-modules.
In particular, it makes sense to consider the integral $\int_Y 1\in\calK_{\GG}$ of the unit cohomology class.

Let us also set $\IH^*_{\GG}(Y)_{loc}=\IH^*_{\GG}(Y)\underset{\calA_{\GG}}\ten \calK_{\GG}$.
Then $\IH^*_{\GG}(Y)_{loc}$ is a finite-dimensional vector space over $\calK_{\GG}$ endowed with a
non-degenerate $\calK_{\GG}$-valued (Poincar\'e) pairing $\la\cdot,\cdot\ra_{Y,\GG}$.

In particular, all of the above is applicable to $Y=\QM^{\theta}_{G,P}$ and $\GG=L\x\CC^*$.
In particular, if we let $1_{G,P}^{\theta}$ denote the unit class in the
$L\x \CC^*$-equivariant cohomology of $\QM_{G,P}^{\theta}$,
then we can define
\eq{partition-finite}
       \calZ_{G,P}=\sum\limits_{\theta\in\Lam_{G,P}^{\theta}}
\grq^{\theta} \ \int\limits_{\QM_{G,P}^{\theta}}1_{G,P}^{\theta}.
\end{equation}
This is a formal series in $\grq\in Z(\chL)$ with values in the
field $\calK_{L\x\CC^*}$ of $L$-invariant rational functions on
$\grl\x\CC$.

In fact the function $\calZ_{G,P}$ is a familiar object in Gromov-Witten theory:
it is explained in \cite{B} that up to a simple factor $\calZ_{G,P}$
is the so called {\it equivariant $J$-function} of $\calG_{G,P}$
(cf. Section 6 of \cite{B}). In particular, this function
was studied from many different points of view ( cf. \cite{GK}, \cite{Kim}, \cite{B} for the case when $P$ is the Borel subgroup).
It was conjectured in \cite{B},\cite{Bicm}
that the function $\calZ_{G,P}$ should have an interpretation in terms of representation theory related
to the Langlands dual Lie algebra $\check\grg$. More generally, let us set
$$
\begin{aligned}
\IH_{G,P}^{\theta}=\IH_{L\x\CC^*}(\QM_{G,P}^{\theta})_{loc},\quad \IH_{G,P}=
\bigoplus\limits_{\theta\in\Lam_{G,P}^+}\IH_{G,P}^{\theta}, \\
\la\cdot,\cdot\ra_{G,P}^{\theta}=\la\cdot,\cdot\ra_{\QM_{G,P}^{\theta},L\x\CC^*},\quad
\la\cdot,\cdot\ra_{G,P}=\bigoplus\limits_{\theta\in\Lam_{G,P}^+}
(-1)^{\langle\theta,\check{\rho}\rangle}\la\cdot,\cdot\ra_{G,P}^{\theta},
\end{aligned}
$$
where $\check{\rho}$ denotes the half-sum of the positive roots of $\grg$,
and we view $\theta$ as the positive integral linear combination of (the
images of) the simple coroots out of ${\mathfrak l}$.

Then one would like to
interpret the $\Lam_{G,P}^+$-graded vector space $\IH_{G,P}$ together with the intersection pairing $\la\cdot,\cdot\ra$ and
the unit cohomology vectors $1_{G,P}^{\theta}\in\IH_{G,P}^{\theta}$ in such terms.
A complete answer to this problem in the case when $P$ is a Borel subgroup was given in \cite{B}.
The purpose of this paper is to suggest a conjectural answer in the general case and to prove this
conjecture for $G=GL(N)$. The relevant representation theory turns out to be the representation theory
of {\em finite $W$-algebras} which we recall in \refs{W-algebras}. We should also note that our conjecture
is motivated by the so called Alday-Gaiotto-Tachikawa (or AGT) conjecture \cite{AGT} which relates 4-dimensional super-symmetric
gauge theory to certain 2-dimensional conformal field theory (more precisely, one may view some part of the
AGT conjecture as an affine version of our conjecture; this point of view is explained in \refs{AGT}).
In fact, the current work grew out of an attempt to create an approach to the AGT conjecture in terms
of geometric representation theory. We hope to pursue this point of view in future publications.
%---------------------------------------------------------
\ssec{main-int}{The main conjecture}In the remainder of this introduction we are going to give a more
 precise formulation of our main conjecture  and indicate the idea of the proof for $G=GL(N)$.
To do this, let us first recall the corresponding result from \cite{B} dealing with the
case when $P$ is a Borel subgroup. In what follows we shall denote it by $B$ instead of $P$.

First, it is shown in \cite{B} that the Lie algebra $\chg$
acts naturally on $\IH_{G,B}$. Moreover, this action has the following
properties. First of all, let us denote by
$\la\cdot,\cdot\ra_{G,B}$ the direct sum of the pairings
$(-1)^{\la \theta,\check{\rho}\ra}\la\cdot,\cdot\ra_{G,B}^{\theta}$.
%\footnote{here $\check{\rho}$ denotes the half-sum of the positive roots
%of $\grg$.}

Recall that the Lie algebra $\chg$ has its triangular
decomposition
$\chg=\chn_+\oplus\chh\oplus\chn_-$. Let $\kap:\chg\to\chg$ denote the
Cartan anti-involution which interchanges $\chn_+$ and $\chn_-$ and
acts as identity on $\chh$.
For each $\lam\in\grh=(\chh)^*$ we denote by $M(\lam)$ the
corresponding
Verma module with highest weight $\lam$; this is a module generated by
a vector $v_\lam$ with (the only) relations
$$
t(v_\lam)=\lam(t)v_\lam\quad\text{for $t\in\chh$ and}\quad
n(v_\lam)=0\quad\text{for $n\in\chn_+$}.
$$

\th{flag}
\begin{enumerate}
\item
 $\IH_{G,B}$ (with the
above action) becomes isomorphic to $M(\lam)$, where
$\lam=-\frac{a}{\hbar}-\rho$.

\item
 $\IH_{G,B}^{\theta}\subset \IH_{G,B}$ is the
$-\frac{a}{\hbar}-\rho-\theta$-weight space of $\IH_{G,B}$.

\item For each $g\in\chg$ and $v,w\in \IH_{G,B}$ we have
$$
\la g(v),w\ra_{G,B}=\la v,\kap(g)w\ra_{G,B}.
$$

\item
 The  vector
$\sum_{\theta}1_{G,B}^{\theta}$ (lying is some completion of
$\IH_{G,B}$) is a Whittaker vector (i.e. a $\grn_+$-eigen-vector) for the above action.
\end{enumerate}
\eth

As a corollary we get that the function
$\grq^{\frac{a}{\hbar}}\calZ_{G,B}$ is an eigen-function of
the {\it quantum Toda hamiltonians}
 associated with $\chg$
with eigen-values determined (in the natural way) by $a$
 (we refer the reader to
\cite{etingof} for the definition of (affine) Toda integrable system
and its relation with Whittaker functions).
In fact, in \cite{B} a similar statement is proved also when $G$ is replaced by the corresponding affine
Kac-Moody group --- cf. \refs{AGT} for more detail.

Our main conjecture gives a generalization of the statements 1)--3) above to arbitrary $P$. Namely,
to any nilpotent element $e\in \chg$ one can associate the so called {\em finite
$W$-algebra} $U(\chg,e)$. We recall the definition in \refss{alg} (this definition is such that
when $e=0$ we have $U(\chg,e)=U\chg$ and when $e$ is regular, then $U(\chg,e)$ is the center of $U\chg$).
Roughly speaking, we conjecture that analogs of 1)--3) hold when $U\chg$ is replaced by $U(\chg,e_{\chL})$ (we refer the reader
to \refs{W-algebras} for the definition of Verma module and Whittaker vectors for finite $W$-algebras). The main purpose of this
paper is to formulate this conjecture more precisely and to prove it for $G$ of type $A$.
%---------------------------------------------------------------------------------------------------
%\ssec{quantum cohomology}{Relation to quantum cohomology of partial flag manifolds and parabolic Toda chains}
%---------------------------------------------------------------------------------------------------
\ssec{organization}{Organization of the paper}The paper is organized as follows:
in \refs{W-algebras} we recall basic definitions about finite $W$-algebras for
general $G$;  we also recall the basic results about parabolic quasi-maps'
spaces and formulate our main conjecture.
In \refs{YW} we recall the results of Brundan and Kleshchev who interpret finite $W$-algebras
in type $A$ using certain shifted Yangians and in \refs{whit} we discuss the notion of Whittaker vectors for finite
$W$-algebras from this point of view. In \refs{parlam} we use it in order to prove our main conjecture for $G=SL(N)$
(and any parabolic). One important ingredient in the proof is this: we replace the intersection cohomology of
parabolic quasi-maps' spaces by the ordinary cohomology of a small resolution of those spaces (which we
call {\em parabolic Laumon spaces}). Finally in \refs{AGT} we discuss the relation between the above results and the AGT
conjecture.
%-----------------------------------------------------------------------------------------
\ssec{}{Acknowledgements} We are grateful to A.~Molev for the explanation of
the results of~\cite{FMO}. We are also grateful to J.~Brundan, A.~Kleschev and I.~Losev
for their explanations about $W$-algebras. Thanks are due to
A.~Tsymbaliuk for the careful reading of the first draft of this note and
spotting several mistakes, and to S.~Gukov, D.~Maulik, A.~Nietzke, A.~Okounkov,
V.~Pestun, Y.~Tachikawa for very helpful discussions on the subject.
A.~B. was partially supported by the NSF grants DMS-0854760 and DMS-0901274.
B.~F, M.~F., and L.~R. were partially
supported by the RFBR grant 09-01-00242, the Ministry of Education and Science of Russian Federation grant No. 2010-1.3.1-111-017-029, 
and the AG Laboratory HSE, RF government grant, ag. 11.G34.31.0023, and
HSE science foundation grant 10-01-0078.

%%%%%%%%%%%%%%%%%%%%%%%%%%%%%%%%%%%%%%%%%%%%%%%%%%%%%%%%%%%%%%%%%%%%%%%%%%%%%%%%%%%%%%%
\sec{W-algebras}{Finite $W$-algebras}

\ssec{alg}{$W$-algebra}
 Let $e$ be a principal
nilpotent element of the Levi subalgebra $\fl\subset\fg$. Let
$U(\fg,e)$ denote the finite $W$-algebra associated to $e$, see
e.g.~\cite{BGK}. We recall its definition for the readers' convenience.
Choose an $\mathfrak{sl}_2$-triple $(e,h,f)$ in $\fg$. We introduce a grading
on $\fg$ by eigenvalues of $\on{ad}_h:$
$$
\fg=\bigoplus_{i\in\BZ}\fg(i),\
\fg(i):=\{\xi\in\fg:\ [h,\xi]=i\xi\}.
$$
 The Killing form $(\cdot,\cdot)$ on
$\fg$ allows to identify $\fg$ with the dual space $\fg^*$. Let $\chi=(e,?)$ be an element of
$\fg^*$ corresponding to $e$. Note that $\chi$ defines a symplectic form
$\omega_\chi$ on $\fg(-1)$ as follows: $\omega_\chi(\xi,\eta):=\langle\chi,
[\xi,\eta]\rangle$. We fix a Lagrangian subspace $l\subset\fg(-1)$ with
respect to $\omega_\chi$, and set $\fm:=l\oplus\bigoplus_{i\leq-2}\fg(i)$.
We define the affine subspace $\fm_\chi\subset U\fg$ as follows:
$\fm_\chi:=\{\xi-\langle\chi,\xi\rangle,\ \xi\in\fm\}$. Finally, we define
the $W$-algebra $U(\fg,e):=(U\fg/U\fg\cdot\fm_\chi)^{\on{ad}\fm}:=
\{a+U\fg\cdot\fm_\chi:\ [\fm,a]\subset U\fg\cdot\fm_\chi\}$.
It is easy to see that
$$
U(\fg,e)=\End_{U\fg}(U\fg\underset{U\grm}\ten \CC_{\chi}),
$$
where $\CC_{\chi}$ denotes the natural 1-dimensional module over $\grm$ corresponding to the character $\chi$.
In this description the algebra structure on $U(\grg,e)$ becomes manifest.

It is equipped with the {\em Kazhdan filtration}
$\on{F}_0U(\fg,e)\subset \on{F}_1U(\fg,e)\subset\ldots$, see
e.g.~section~3.2 of~\cite{BGK}. We recall its definition for the readers'
convenience. We denote the standard PBW filtration on $U\fg$ (by the order
of a monomial) by $\on{F}_i^{st}U\fg$. The Kazhdan filtration on $U\fg$
is defined by $\on{F}_iU\fg:=\sum_{2k+j\leq i}\on{F}_k^{st}U\fg\cap U\fg(j)$
where $U\fg(j)$ is the eigenspace of $\on{ad}_h$ on $U\fg$ with eigenvalue $j$.
Being a subquotient of $U\fg$, the $W$-algebra $U(\fg,e)$ inherits the Kazhdan
filtration $\on{F}_iU(\fg,e)$.

We also consider the {\em shifted
Kazhdan filtration} $F_iU(\fg,e):=\on{F}_{i+1}U(\fg,e)$, and we
define the $\BC[\hbar]$-algebra $U^\hbar(\fg,e)$ as the {\em Rees
algebra} of the filtered algebra $(U(\fg,e),F_\bullet)$. Abusing
notation, we will sometimes call $U^\hbar(\fg,e)$ just a
$W$-algebra.

%Let us choose a good grading of $\fg$ associated with $e$ (see~\cite{BGK}).
%Let $\fl'\subset\fp'$ be the zero graded part, and the nonnegative graded
%part, respectively. It is a Levi and a parabolic subalgebra of $\fg$.
Let us extend the scalars to the field $\CK:=\CK_{\check{L}\times\BC^*}$, and let
$\hbar\in\CK$ stand for the generator of $H^2_{\BC^*}(pt,\BZ)$.
Thus $\CK=\BC(\ft^*/W_L\times\BA^1)$ where $W_L$ stands for the Weyl group of
$\grl$, and $\BA^1$ is the affine line with coordinate $\hbar$.

Given a point $\Lambda\in\ft^*/W_L$ we consider the Verma module
$M(-\hbar^{-1}\Lambda,e)$ over $U^\hbar(\fg,e)\otimes\CK\simeq
U(\fg,e)\otimes\CK$ introduced in sections~4.2 and~5.1 of~\cite{BGK}.
As $\Lambda\in\ft^*/W_L$ varies, these modules form a family over
$\BC[\ft^*/W_L\times(\BA^1-0)]$. Localizing to $\CK$ we obtain the
{\em universal Verma module} $M(\fg,e)$ over $U(\fg,e)\otimes\CK$.
Note that a certain $\rho$-shift is incorporated into the definition of
$M(\fg,e)$, cf. the text right after~Lemma~5.1 of~\cite{BGK}
\footnote{The definition of the module $M(\fg,e)$ from \cite{BGK} is
(unfortunately) quite involved and
we are not going to recall it here. On the other hand for $\grg=\fgl(N)$ there is another (in some sense, more explicit)
definition of this module which we are going to recall in \refs{YW}.}.

In what follows we shall often abbreviate $M:=M(\fg,e)$.
%-------------------------------------------------------------------------------------------------
\ssec{wh}{Whittaker vector}
Let $\ft^e$ stand for the centralizer of $e$ in $\ft$. Recall from~\cite{BGK}
that the collection of nonzero weights of $\ft^e$ on $\fg$ is called a
{\em restricted root system} $\Phi^e$, and the weights on
$\fp':=\bigoplus_{i\geq0}\fg(i)$ form a
positive root system $\Phi^e_+\subset\Phi^e$. Let $I^e$ be the set of simple
roots, i.e. positive roots which are not positive linear combinations of other
positive roots. According to~Theorem~6 of~\cite{BG}, the simple roots form
a base of $(\ft^e)^\vee$.

Let us choose a linear embedding $\Theta:\ \fg^e\hookrightarrow
U^\hbar(\fg,e)$ as in~Theorem~3.6 of~\cite{BGK}. For a simple root
$\alpha\in I^e$, we consider the corresponding weight space
$\Theta(\fg^e_\alpha)$. The Kazhdan filtration induces the increasing
filtration on the root space $\Theta(\fg^e_\alpha)$. We define a positive
integer $m_\alpha$ so that $\on{F}_{m_\alpha}\Theta(\fg^e_\alpha)=
\Theta(\fg^e_\alpha)$, but $\on{F}_{m_\alpha-1}\Theta(\fg^e_\alpha)\ne
\Theta(\fg^e_\alpha)$. The following conjecture holds for $\fg$ of type $A$
by the work of J.~Brundan and A.~Kleshchev (cf.~\refss{positive} below),
and for all
exceptional types according to computer calculations by J.~Brundan (private
communication):

\conj{multip}
$\dim\on{F}_{m_\alpha}\Theta(\fg^e_\alpha)/\on{F}_{m_\alpha-1}\Theta(\fg^e_\alpha)=1$.
\econj

\defe{regu}
A linear functional $\psi$ on $\bigoplus_{\alpha\in I^e}\Theta(\fg^e_\alpha)$
is called regular if for any $\alpha\in I^e$ we have
$\psi(\on{F}_{m_\alpha-1}\Theta(\fg^e_\alpha))=0$ but
$\psi(\on{F}_{m_\alpha}\Theta(\fg^e_\alpha))\ne0$.
\edefe

If~\refco{multip} is true, then $T^e$ acts simply transitively on the set of
regular functionals.

\defe{regul}
Given a regular functional $\psi$ on
$\bigoplus_{\alpha\in I^e}\Theta(\fg^e_\alpha)$, a $\psi$-eigenvector $\fw$
in a completion $\prod_{\theta\in\Lambda_{{\check G},{\check P}}}M_\theta$ of
the universal Verma module $M$ is called a $\psi$-Whittaker vector.
\edefe

\ssec{form}{Shapovalov form}
Let $\sigma$ stand for the Cartan antiinvolution of $\fg$ identical on $\ft$.
Let $w_0^\fl$ stand for the adjoint action of a representative of the longest
element of the Weyl group of the Levi subalgebra $\fl$. Then the composition
$w_0^\fl\sigma$ preserves $e$ and everything else entering the definition of
the finite $W$-algebra and gives rise to an antiisomorphism
$U^\hbar(\fg,e)\iso\overline{U}{}^\hbar(\fg,e)$ where
$\overline{U}{}^\hbar(\fg,e)$ (see~Section~2.2 of~\cite{BGK}) is defined just as
$U^\hbar(\fg,e)$, only with left ideals replaced by right ideals. Composing
this antiisomorphism with the isomorphism $\overline{U}{}^\hbar(\fg,e)\iso
U^\hbar(\fg,e)$ of~Corollary~2.9 of~\cite{BGK} we obtain an antiinvolution
$\varsigma$ of $U^\hbar(\fg,e)$.

\defe{for} The Shapovalov bilinear form $(\cdot,\cdot)$ on the universal
Verma module $M$ with values in $\CK$ is the unique
bilinear form such that $(x,yu)=(\varsigma(y)x,u)$ for any $x,u\in M,\
y\in U^\hbar(\fg,e)$ with value 1 on the highest vector.
\edefe

%%%%%%%%%%%%%%%%%%%%%%%%%%%%%%%%%%%%%%%%%%%%%%%%%%%%%%%%%%%%%%%%%%%%%%%%%%%%%%%%%%%%%%%%%%%%%%%%%%%%%%
\ssec{main-general}{Parabolic Zastava spaces and finite $W$-algebras: the
main conjecture}
We are now ready to formulate our main conjectures. We are going to change slightly the notations. Namely,
the symbols $G,P,L,\grg,\grp,\grl$ will
denote the same things as in the introduction. However, we are now going to denote by $e$ the principal
nilpotent element in the {\em Langlands dual} Lie algebra
${\check\grl}\subset\grg$; similarly, $M$ will now denote the universal
Verma module over $U(\chg,e)$. Note that in this case $M$ becomes naturally graded by
$\Lam_{G,P}^+$. Also, the universal coefficient field $\CK$ is now nothing else but the field $\CK_{L\x\CC^*}$
which appeared in the introduction. With these conventions we formulate the following
%-----------------------------------------------------------------------------------------------
\conj{main}
\begin{enumerate}
\item
There is an isomorphism $\Psi$ of
$\Lambda_{G,P}$-graded $\CK$-vector spaces $\IH_{G,P}$ and
$M$; in particular, there is a natural $U^{\hbar}(\chg,e)$-action on $\IH_{G,P}$.
\item
The isomorphism $\Psi$ takes the vector $\sum_{\theta\in\Lambda_{G,P}}1^\theta\in
\prod_{\theta\in\Lambda_{G,P}}\IH^\theta_{G,P}$ to a $\psi$-Whittaker vector
$\fw\in\prod_{\theta\in\Lambda_{G,P}}M_\theta$
for certain regular functional $\psi$.
\item
For $x,u\in
V_\theta$ we have $(\Psi(x),\Psi(u))=(-1)^{|\theta|}\langle
x,u\rangle^\theta_{G,P}$.
\end{enumerate}
\econj

In what follows we shall give a more precise formulation of this conjecture when $\grg$ is of type
$A$ (in that case we shall also prove the conjecture). In particular, we shall specify the
regular functional $\psi$ in that case.
%%%%%%%%%%%%%%%%%%%%%%%%%%%%%%%%%%%%%%%%%%%%%%%%%%%%%%%%%%%%%%%%%%%%%%%%%%%%%%%%%%%%%%%%%

%%%%%%%%%%%%%%%%%%%%%%%%%%%%%%%%%%%%%%%%%%%%%%%%%%%%%%%%%%%%%%%%%%%%%%%%%%%%%%%%%%%%%%%%
\sec{YW}{Shifted Yangians and finite $W$-algebras}
In this Section we recall an explicit realization of finite $W$-algebras
in type $A$ using shifted Yangians (due to Brundan and Kleshchev).

\ssec{relations}{Shifted Yangian} Let $\pi=(p_1,\ldots,p_n)$ where
$p_1\leq p_2\leq\ldots\leq p_n$, and $p_1+\ldots+p_n=N$. Recall
the shifted Yangian $Y^\hbar_\pi(\fgl_n)$ introduced by J.~Brundan
and A.~Kleshchev (see~\cite{BK} and~\cite{FMO}). It is an
associative $\BC[\hbar]$-algebra with generators $\sd_i^{(r)},\
i=1,\ldots,n,\ r\geq1;\ \ssf_i^{(r)},\ i=1,\ldots,n-1,\ r\geq1;\
\se_i^{(r)},\ i=1,\ldots,n-1,\ r\geq p_{i+1}-p_i+1,$ subject to the
following relations

\eq{a} [\sd_i^{(r)},\sd_j^{(s)}]=0,
\end{equation}

\eq{b}
[\se_i^{(r)},\ssf_j^{(s)}]=-\delta_{ij}\hbar\sum_{t=0}^{r+s-1}\sd'_i{}^{(t)}
\sd_{i+1}^{(r+s-t-1)},
\end{equation}
where $\sd_i^{(0)}:=1$, and the elements $\sd'_i{}^{(r)}$ are found
from the relations
$\sum_{t=0}^r\sd_i^{(t)}\sd'_i{}^{(r-t)}=\delta_{r0},\
r=0,1,\ldots$;

\eq{c} [\sd_i^{(r)},\se_j^{(s)}]=
\hbar(\delta_{ij}-\delta_{i,j+1})\sum_{t=0}^{r-1}\sd_i^{(t)}\se_j^{(r+s-t-1)},
\end{equation}

\eq{d}
[\sd_i^{(r)},\ssf_j^{(s)}]=\hbar(\delta_{i,j+1}-\delta_{ij})\sum_{t=0}^{r-1}
\ssf_j^{(r+s-t-1)}\sd_i^{(t)},
\end{equation}

\eq{e} [\se_i^{(r)},\se_i^{(s+1)}]-[\se_i^{(r+1)},\se_i^{(s)}]=
\hbar(\se_i^{(r)}\se_i^{(s)}+\se_i^{(s)}\se_i^{(r)}),
\end{equation}

\eq{f} [\ssf_i^{(r+1)},\ssf_i^{(s)}]-[\ssf_i^{(r)},\ssf_i^{(s+1)}]=
\hbar(\ssf_i^{(r)}\ssf_i^{(s)}+\ssf_i^{(s)}\ssf_i^{(r)}),
\end{equation}

\eq{g}
[\se_i^{(r)},\se_{i+1}^{(s+1)}]-[\se_i^{(r+1)},\se_{i+1}^{(s)}]=
-\hbar\se_i^{(r)}\se_{i+1}^{(s)},
\end{equation}

\eq{h}
[\ssf_i^{(r+1)},\ssf_{i+1}^{(s)}]-[\ssf_i^{(r)},\ssf_{i+1}^{(s+1)}]=
-\hbar\ssf_{i+1}^{(s)}\ssf_i^{(r)},
\end{equation}

\eq{i} [\se_i^{(r)},\se_j^{(s)}]=0\qquad\on{if}\ |i-j|>1,
\end{equation}

\eq{j} [\ssf_i^{(r)},\ssf_j^{(s)}]=0\qquad\on{if}\ |i-j|>1,
\end{equation}

\eq{k}
[\se_i^{(r)},[\se_i^{(s)},\se_j^{(t)}]]+[\se_i^{(s)},[\se_i^{(r)},\se_j^{(t)}]]
=0\qquad\on{if}\ |i-j|=1,
\end{equation}

\eq{l} [\ssf_i^{(r)},[\ssf_i^{(s)},\ssf_j^{(t)}]]+
[\ssf_i^{(s)},[\ssf_i^{(r)},\ssf_j^{(t)}]]=0\qquad\on{if}\ |i-j|=1.
\end{equation}

We introduce the generating series
$$
\begin{aligned}
\sd_k(u)=1+\sum_{s=1}^\infty\sd_k^{(s)}\hbar^{-s+1}u^{-s}, &\quad
\se_k(u)=\sum_{s=p_{k+1}-p_k+1}^\infty\se_k^{(s)}\hbar^{-s+1}u^{-s},\\
\ssf_k(u):=&\sum_{s=1}^\infty\ssf_k^{(s)}\hbar^{-s+1}u^{-s}.
\end{aligned}
$$
Finally, we define $\sa_k(u):=\sd_1(u)\sd_2(u-1)\ldots\sd_k(u-k+1),\
\sA_k(u):=u^{p_1}(u-1)^{p_2}\ldots(u-k+1)^{p_k}\sa_k(u)$, and also

\eq{dvass}
\sB_k(u):=(u-k+1)^{p_{k+1}-p_k}\sA_k(u)\se_k(u-k+1),
\end{equation}

\eq{triss} \sC_k(u):=\ssf_k(u-k+1)\sA_k(u).
\end{equation}

\ssec{w-algebra}{$W$-algebra $U^\hbar(\fg,e)$ and its universal
Verma module} Let $\fp\subset\fg:=\fgl_N=\End(W)$ be the parabolic
subalgebra preserving the flag $0\subset W_1\subset\ldots\subset
W_{n-1}\subset W$. Here $W=\langle w_1,\ldots,w_N\rangle$, and
$W_i=\langle w_1,\ldots,w_{p_1+\ldots+p_i}\rangle$.
Let $e$ be a principal nilpotent element of a
Levi factor $\fl$ of $\fp$. Let $U(\fg,e)$ denote the finite
$W$-algebra associated to $e$, see e.g.~\cite{BK}. It is equipped
with the {\em Kazhdan filtration} $\on{F}_0U(\fg,e)\subset
\on{F}_1U(\fg,e)\subset\ldots$, see e.g.~section~3.2 of~\cite{BK}.
We also consider the {\em shifted Kazhdan filtration}
$F_iU(\fg,e):=\on{F}_{i+1}U(\fg,e)$, and we define the
$\BC[\hbar]$-algebra $U^\hbar(\fg,e)$ as the {\em Rees algebra} of
the filtered algebra $(U(\fg,e),F_\bullet)$. Abusing notation, we
will sometimes call $U^\hbar(\fg,e)$ just a $W$-algebra.

According to~section~3.4 of~\cite{BK}, 
$U^\hbar(\fg,e)\otimes_{\BC[\hbar]}\BC(\hbar)$ is the
quotient of $Y^\hbar_\pi(\fgl_n)\otimes_{\BC[\hbar]}\BC(\hbar)$ 
by the relations $\sd_1^{(r)}=0,\ r>p_1$.

%Moreover, according to~section~3.6 of~\cite{BK}, there is an
%injective homomorphism ({\em Miura transform})
%$U^\hbar(\fg,e)\hookrightarrow U(\fl')$ where $\fl'\subset\fgl_N$ is
%another Levi subalgebra, the zero grading part of a good grading of
%$\fg$ associated to $e$. It is a Levi factor of a parabolic
%subalgebra $\fq$: the nonnegative grading part of a good grading of
%$\fg$ associated to $e$. We have $\fl'\simeq\fgl_{q_1}\oplus\ldots
%\oplus\fgl_{q_l}$ where $(q_1\geq q_2\geq\ldots\geq q_l)$ is the
%partition of $N$ dual to $\pi$. We define the {\em universal Verma
%module} $M$ over $U^\hbar(\fg,e)$ as the restriction of the
%universal Verma module over $\fl'$ to $U^\hbar(\fg,e)$.
Let us denote the standard coordinates in the diagonal Cartan subalgebra $\ft$
of $\fgl_N$ by $x_1,\ldots,x_N$. Then the universal Verma module
$M=M(\fgl_N,e)$ over
$U^\hbar(\fg,e)\otimes\CK$ is a vector space over the field $\CK$ of
rational functions in $\hbar,x_1,\ldots,x_N$ symmetric in the groups
$(x_1,\ldots,x_{p_1}),(x_{p_1+1},\ldots,x_{p_1+p_2}),\ldots,(x_{p_1+\ldots+p_{n-1}+1},
\ldots,x_N)$.
Let us consider the field extension $\CK\subset\CK':=\BC(\ft^*\times\BA^1)=
\BC(\hbar,x_1,\ldots,x_N)$. We denote by $M'=M'(\fgl_N,e):=M\otimes_\CK\CK'$
the universal Verma module with the extended scalars.

\ssec{Gelfand}{The Gelfand-Tsetlin module}
According to~\cite{FMO}, the module $M'$ admits a rather explicit description.
More precisely, the authors of~\cite{FMO} construct a
{\em Gelfand-Tsetlin module} $V$ over $U^\hbar(\fg,e)\otimes\CK'$ equipped
with a {\em Gelfand-Tsetlin base}
numbered by the {\em Gelfand-Tsetlin patterns}, and write down
explicitly the matrix coefficients of the generators of
$U^\hbar(\fg,e)$ in this base. We recall these results here. To a
collection $\vec{\unl{d}}=(d_{ij}^{(a)}),\ n-1\geq i\geq j,\ p_j\geq
a\geq1$, we associate a Gelfand-Tsetlin pattern
$\Lambda=\Lambda(\vec{\unl{d}}):=(\lambda_{ij}^{(a)}),\ n\geq i\geq
j,\ p_j\geq a\geq1$, as follows:
$\lambda_{nj}^{(a)}:=\hbar^{-1}x_{p_1+\ldots+p_{j-1}+a}+j-1,\ n\geq
j\geq 1;\ \lambda_{ij}^{(a)}:=-\hbar^{-1}\sfp_{ij}^{(a)}+j-1,\
n-1\geq i\geq j\geq1$ (see~\refe{pij}). The corresponding base
element $\xi_\Lambda= \xi_{\Lambda(\vec{\unl{d}})}$ will be denoted
by $\xi_{\vec{\unl{d}}}$ for short. Thus, the set
$\{\xi_{\vec{\unl{d}}}\}$ (over all collections $\vec{\unl{d}}$)
forms a basis of $V$.

We have used the following notation: given $\vec{\unl{d}}$,
\eq{pij}
\sfp_{ik}^{(l)}:=\hbar d_{ik}^{(l)}-x_{p_1+\ldots+p_{k-1}+l},\ 1\leq l\leq p_i
\end{equation}
Also, for $n\geq i\geq j$ we introduce the monic polynomials
$\lambda_{ij}(u):=
(u+\lambda_{ij}^{(1)})\ldots(u+\lambda_{ij}^{(p_j)})$.

Finally, we define the action of the generators of $U^\hbar(\fgl_N,e)$
on $V$ by their matrix elements in the Gelfand-Tsetlin base:
\eq{eee}
\se_{i[\vec{\unl{d}},\vec{\unl{d}}{}']}^{(s)}=
-\hbar^{-1-p_i}(\sfp_{ij}^{(a)}-i\hbar)^{s-1-p_{i+1}+p_i}
\prod_{\substack{k\leq i,\ b\leq p_k\\ (k,b)\ne(j,a)}}
(\sfp_{ij}^{(a)}-\sfp_{ik}^{(b)})^{-1}
\prod_{k\leq i+1}\prod_{b\leq p_k}(\sfp_{ij}^{(a)}-\sfp_{i+1,k}^{(b)}),
\end{equation}
if $d_{ij}^{(a)}{}'=d_{ij}^{(a)}-1$ for certain $j\leq i$;
\eq{fff}
\ssf_{i[\vec{\unl{d}},\vec{\unl{d}}{}']}^{(s)}=
\hbar^{-1+p_i}(\sfp_{ij}^{(a)}+(1-i)\hbar)^{s-1}
\prod_{\substack{k\leq i,\ b\leq p_k\\ (k,b)\ne(j,a)}}
(\sfp_{ij}^{(a)}-\sfp_{ik}^{(b)})^{-1}
\prod_{k\leq i-1}\prod_{b\leq p_k}(\sfp_{ij}^{(a)}-\sfp_{i-1,k}^{(b)}),
\end{equation}
if $d_{ij}^{(a)}{}'=d_{ij}^{(a)}+1$ for certain $j\leq i$.
All the other matrix coefficients of $\se_i^{(s)},\ssf_i^{(s)}$ vanish.

The following Proposition is taken from~\cite{FMO}.
\prop{fmo}
The formulas~\refe{eee},~\refe{fff} give rise to the action of
$U^\hbar(\fgl_N,e)\otimes\CK'$
on $V$. Moreover, $\sA_i(u), \sC_i(u), \sB_i(u)$ are polynomials in $u$
of degrees $p_1+\ldots+p_i$ (resp. $p_1+\ldots+p_i-1$, $p_1+\ldots+p_i-1$)
such that
\eq{A}
\sA_i(u)\xi_{\vec{\unl{d}}}=\lambda_{i1}(u)\ldots\lambda_{ii}(u-i+1)
\xi_{\vec{\unl{d}}},
\end{equation}
for $i=1,\ldots,n$, and
\eq{B}
\sB_i(\hbar^{-1}\sfp_{ij}^{(a)})\xi_{\vec{\unl{d}}}=-\lambda_{i+1,1}
(\hbar^{-1}\sfp_{ij}^{(a)})\lambda_{i+1,2}(\hbar^{-1}\sfp_{ij}^{(a)}-1)\ldots
\lambda_{i+1,i+1}(\hbar^{-1}\sfp_{ij}^{(a)}-i)\xi_{\vec{\unl{d}}+\delta_{ij}^{(a)}},
\end{equation}

\eq{C}
\sC_i(\hbar^{-1}\sfp_{ij}^{(a)})\xi_{\vec{\unl{d}}}=\lambda_{i-1,1}
(\hbar^{-1}\sfp_{ij}^{(a)})\lambda_{i-1,2}(\hbar^{-1}\sfp_{ij}^{(a)}-1)\ldots
\lambda_{i-1,i-1}(\hbar^{-1}\sfp_{ij}^{(a)}-i+2)\xi_{\vec{\unl{d}}-\delta_{ij}^{(a)}},
\end{equation}
for $i=1,\ldots,n-1$, where $\vec{\unl{d}}\pm\delta_{ij}^{(a)}$ is
obtained from $\vec{\unl{d}}$ by replacing $d_{ij}^{(a)}$ by
$d_{ij}^{(a)}\pm1$.
\eprop

\prf
The formulas~\refe{B} and~\refe{eee} (resp.~\refe{C} and~\refe{fff})
are equivalent by the Lagrange interpolation. So if suffices to prove that
the formulas~\refe{A},~\refe{B},~\refe{C} give rise to the action of
$U^\hbar(\fgl_N,e)\otimes\CK'$ on $V$. Now $V$ admits a certain integral form
over $\BC[\hbar^{\pm1},x_1,\ldots,x_N]$ which can be specialized to the values
of parameters $x_i$ satisfying certain integrality and positivity conditions.
These specializations admit finite-dimensional subspaces spanned by certain
finite subsets of the Gelfand-Tsetlin base. Theorem~4.1 of~\cite{FMO} describes
the action of $U^\hbar(\fgl_N,e)$ in these finite-dimensional subspaces by
the formulas~\refe{A},~\refe{B},~\refe{C}. It follows that the action of
generators given by~\refe{A},~\refe{B},~\refe{C} in these subspaces satisfies
the relations of $U^\hbar(\fgl_N,e)$. For each given $\vec{\unl{d}}$, both
$\vec{\unl{d}}$ and $\vec{\unl{d}}\pm\delta_{ij}^{(a)}$ enter the above
finite subsets for quite a few specializations: more precisely, the set of
special values of $x_1,\ldots,x_N$ such that both $\vec{\unl{d}}$ and
$\vec{\unl{d}}\pm\delta_{ij}^{(a)}$ enter the corresponding finite subsets
is Zariski dense in $\ft^*\times\BA^1$. It follows that the relations of
$U^\hbar(\fgl_N,e)$ are satisfied for all values of $x_1,\ldots,x_N$.
\epr

\prop{GelTs}
The Gelfand-Tsetlin module $V$ is isomorphic to the universal Verma module
with extended scalars $M'$.
\eprop

\prf
First, $V$ is an irreducible module over $U^\hbar(\fgl_N,e)\otimes\CK'$.
In effect, by~\refp{fmo} the Gelfand-Tsetlin subalgebra of
$U^\hbar(\fgl_N,e)\otimes\CK'$ generated by $\sd^{(r)}_i$ acts diagonally
in the Gelfand-Tsetlin base with pairwise distinct eigenvalues. Therefore,
it suffices to check the following two things:

(1) for each $\vec{\unl{d}}$ there are indices $i,s$ such that
$\ssf_i^{(s)}\xi_{\vec{\unl{d}}}\ne0$;

(2) for each $\vec{\unl{d}}\ne0$ there are indices $i,s$ such that
$\se_i^{(s)}\xi_{\vec{\unl{d}}}\ne0$.

Both follow directly from the formulas~\refe{eee},~\refe{fff}.

Second, $M'$ is an irreducible module over $U^\hbar(\fgl_N,e)\otimes\CK'$.
In effect, for a general highest weight  $\Lambda\in\ft^*/W_L$ the Verma module
$M(-\hbar^{-1}\Lambda,e)$ is irreducible according to~\cite{BGK}. Hence the
universal Verma module is irreducible as well.

Now to construct the desired isomorphism $M'\iso V$ it suffices to produce
a nonzero homomorphism $M'\to V$. By the universal property of Verma modules,
it suffices to identify the highest weights of $M'$ and $V$.
By the argument in the proof of~Theorem~5.5 of~\cite{BGK}, $\sd_i^{(r)}$
acts on the highest vector of $M'$ by multiplication by the $r$th elementary
symmetric polynomial in the variables $i-1-\hbar^{-1}x_j,\ p_1+\ldots+p_{i-1}+1
\leq j\leq p_1+\ldots+p_i$. On the other hand, it follows from the
formula~\refe{A} that $\sd_i(u)$ acts on the highest vector $\xi_0$ by
multiplication by
$u^{-p_i}\prod_{j=p_1+\ldots+p_{i-1}+1}^{p_1+\ldots+p_i}(u+i-1-\hbar^{-1}x_j)$.
The coincidence of highest weights completes the proof of the proposition.
\epr

Recall that the Galois group of $\CK'$ over $\CK$ is $W_L$. By the
irreducibility of $U^\hbar(\fgl_N,e)\otimes\CK'$-module $V$, there is a unique
{\em semilinear} action of $W_L$ on $V$ intertwining the action of
$U^\hbar(\fgl_N,e)\otimes\CK'$, and trivial on the highest vector $\xi_0$.

\cor{VM}
The universal Verma module $M(\fgl_N,e)$ is isomorphic to $V^{W_L}$.
\ecor

\ssec{positive}{The characters of positive subalgebra}
Recall the notations of~\cite{BGK}. So $\ft$ stands for the diagonal Cartan
subalgebra of $\fg$, and $\ft^e$ stands for the centralizer of $e$ in $\ft$,
and $\fg^e$ stands for the centralizer of $e$ in $\fg$. The collection $\Phi^e$
of non-zero weights of $\ft^e$ on $\fg$ is a {\em restricted root system},
see e.g.~section~3.1 of~\cite{BGK}. The roots appearing in $\fp'\subset\fg$
form a system $\Phi^e_+\subset\Phi^e$ of positive roots. Let us denote by
$\fg^e_+\subset\fg^e$ the subspace spanned by the positive root vectors.
Recall a linear space embedding $\Theta:\ \fg^e\hookrightarrow U^\hbar(\fg,e)$
of~Theorem~3.6 of~\cite{BGK}. We define $U^\hbar_+(\fg,e)$ as the subalgebra
of $U^\hbar(\fg,e)$ generated by $\Theta(\fg^e_+)$.
%It is independent of a choice of $\Theta$.
In terms of the shifted Yangian, $U^\hbar_+(\fg,e)$ is
generated by $\{\se_i^{(s)}\},\ 1\leq i\leq n-1,\ s\geq p_{i+1}-p_i+1$.

We are interested in the (additive) characters of $U^\hbar_+(\fg,e)$, that is
maximal ideals of $U^\hbar_+(\fg,e)^{\on{ad}}:=
U^\hbar_+(\fg,e)/[U^\hbar_+(\fg,e),U^\hbar_+(\fg,e)]$. We have
$U^\hbar_+(\fg,e)^{\on{ad}}\simeq\on{Sym}(\Theta(\fg^e_+)^{\on{ad}})$
where $\Theta(\fg^e_+)^{\on{ad}}:=\Theta(\fg^e_+)/\Theta(\fg^e_+)\cap
[U^\hbar_+(\fg,e),U^\hbar_+(\fg,e)]$. In terms of roots,
$\Theta(\fg^e_+)^{\on{ad}}$ is spanned by the {\em simple} positive root
spaces in $\Theta(\fg^e_+)$. In terms of the shifted Yangian,
$\Theta(\fg^e_+)^{\on{ad}}$ is spanned by
$\{\se_i^{(s)}\},\ 1\leq i\leq n-1,\ p_{i+1}-p_i+1\leq s\leq p_{i+1}$.

The Kazhdan filtration induces the increasing filtration on the root space
$\Theta(\fg^e_\alpha)$ for a simple positive root $\alpha\in\Phi^e_+$.
In terms of the shifted Yangian, $\on{F}_r\Theta(\fg^e_{\alpha_i})$ is spanned
by $\{\se_i^{(s)}\},\ 1\leq i\leq n-1,\ r\geq s\geq p_{i+1}-p_i+1$.
For a simple root
$\alpha$ we define $m_\alpha$ so that $\on{F}_{m_\alpha}\Theta(\fg^e_\alpha)=
\Theta(\fg^e_\alpha)$, but $\on{F}_{m_\alpha-1}\Theta(\fg^e_\alpha)\ne
\Theta(\fg^e_\alpha)$. Clearly, for $\alpha=\alpha_i$ we have
$m_{\alpha_i}=p_{i+1}$.

We say that an additive character $\chi:\ U^\hbar_+(\fg,e)\to\BC(\hbar)$,
that is a linear function $\Theta(\fg^e_+)^{\on{ad}}\to\BC(\hbar)$,
is {\em regular} if $\chi(\on{F}_{m_\alpha-1}\Theta(\fg^e_\alpha))=0$, but
$\chi(\on{F}_{m_\alpha}\Theta(\fg^e_\alpha))\ne0$ for any simple root $\alpha$.
Let $T^e$ stand for the centralizer of $e$ in the diagonal torus of $G$.
Then the adjoint action of $T^e$ on the set of regular characters is
transitive. We specify one particular regular character in terms of the
shifted Yangian:
$\chi_\hbar(\se_i^{(s)})=0$ for $p_{i+1}-p_i+1\leq s<p_{i+1}$, and
$\chi_\hbar(\se_i^{(p_{i+1})})=\hbar^{-1}$.

\defe{whitt} The Whittaker vector $\fw\in\widehat M$ in a completion of the
universal Verma module (the product of the weight spaces) is the unique
eigenvector for $U^\hbar_+(\fg,e)$ with the eigenvalue $\chi_\hbar$
whose highest weight component coincides with the highest vector.
For a weight $\unl{d}$ we denote by $\fw_{\unl{d}}$ the weight $\unl{d}$
component of $\fw$.
\edefe
\ssec{}{The Shapovalov form in terms of shifted Yangians}
We consider the antiinvolution $\varsigma:\ Y_\pi^\hbar(\fgl_n)\to
Y_\pi^\hbar(\fgl_n)$ taking $\sd_i^{(s)}$ to $\sd_i^{(s)}$, and
$\ssf_i^{(s)}$ to $\se_i^{(s+p_{i+1}-p_i)}$. This is nothing else than the
composition of isomorphism~(2.35) of~\cite{BK} and anti-isomorphism~(2.39)
of~\cite{BK}. It descends to the same named antiinvolution
$\varsigma:\ U^\hbar(\fg,e)\to U^\hbar(\fg,e)$. According to~section~3.5
of~\cite{BK}, this antiinvolution can be alternatively described as follows.
Let $\sigma$ stand for the Cartan antiinvolution of $\fg$ (transposition).
Let $w_0^\fl$ stand for the adjoint action of a representative of the longest
element of the Weyl group of the Levi subalgebra $\fl$. Then the composition
$w_0^\fl\sigma$ preserves $e$ and everything else entering the definition of
the finite $W$-algebra and gives rise to an antiisomorphism
$U^\hbar(\fg,e)\iso\overline{U}{}^\hbar(\fg,e)$ where
$\overline{U}{}^\hbar(\fg,e)$ (see~Section~2.2 of~\cite{BGK}) is defined just as
$U^\hbar(\fg,e)$, only with left ideals replaced by right ideals. Composing
this antiisomorphism with the isomorphism $\overline{U}{}^\hbar(\fg,e)\iso
U^\hbar(\fg,e)$ of~Corollary~2.9 of~\cite{BGK} we obtain an antiinvolution
of $U^\hbar(\fg,e)$. This antiinvolution coincides with $\varsigma$.

\defe{shap} The Shapovalov bilinear form $(\cdot,\cdot)$ on the universal
Verma module $M$ with values in $\BC(\hbar,x_1,\ldots,x_n)$ is the unique
bilinear form such that $(x,yu)=(\varsigma(y)x,u)$ for any $x,u\in M,\
y\in U^\hbar(\fg,e)$, with value 1 on the highest vector.
\edefe

%%%%%%%%%%%%%%%%%%%%%%%%%%%%%%%%%%%%%%%%%%%%%%%%%%%%%%%%%%%%%%%%%%%%%%%%%%%%%%%%%%%%%%%

\sec{parlam}{Parabolic Laumon spaces and correspondences: proof of the main conjecture for $G=GL(n)$}
%%%%%%%%%%%%%%%%%%%%%%%%%%%%%%%%%%%%%%%%%%%%%%%%%%%%%%%%%%%%%%%%%%%%%%%%%%%%%%%%%%%%%%%%%%%%%%%%%%
In this section we prove~\refco{main}  for $G=GL_N$.
Note that in this case $G\simeq{\check G},\ P\simeq{\check P},\
L\simeq{\check L},\ T\simeq{\check T}$.

\ssec{1}{}
We recall the setup of~\cite{FNR}. Let $\bC$ be a smooth
projective curve of genus zero. We fix a coordinate $z$ on $\bC$,
and consider the action of $\BC^*$ on $\bC$ such that
$v(z)=v^{-1}z$. We have $\bC^{\BC^*}=\{0,\infty\}$.

We consider an $N$-dimensional vector space $W$ with a basis
$w_1,\ldots,w_N$. This defines a Cartan torus $T\subset
G=GL_N=Aut(W)$ acting on $W$ as follows: for
$T\ni\unl{t}=(t_1,\ldots,t_N)$ we have
$\unl{t}(w_i)=t_iw_i$.

\ssec{2}{}
We fix an $n$-tuple of positive integers $p_1\leq p_2\leq\ldots\leq p_n$
such that $p_1+\ldots+p_n=N$. Let $P\subset G$ be a parabolic subgroup
preserving the flag $0\subset W_1:=\langle w_1,\ldots,w_{p_1}\rangle\subset
W_2:=\langle w_1,\ldots,w_{p_1+p_2}\rangle\subset\ldots\subset W_{n-1}:=
\langle w_1,\ldots w_{p_1+\ldots+p_{n-1}}\rangle\subset W_n:=W$.
Let $G/P$ be the corresponding partial flag variety.

Given an $(n-1)$-tuple of nonnegative integers
$\unl{d}=(d_1,\ldots,d_{n-1})$, we consider the Laumon's parabolic
quasiflags' space $\CQ_{\unl{d}}$, see~\cite{L1},~4.2. It is the
moduli space of flags of locally free subsheaves
$$0\subset\CW_1\subset\ldots\subset\CW_{n-1}\subset\CW=W\otimes\CO_\bC$$
such that $\on{rank}(\CW_k)=p_1+\ldots+p_k$, and $\deg(\CW_k)=-d_k$.

It is known to be a smooth connected projective variety
of dimension $\dim(G/P)+\sum_{i=1}^{n-1}d_i(p_i+p_{i+1})$, see~\cite{L1},~2.10.

\ssec{3}{} We consider the following locally closed subvariety
$\fQ_{\unl{d}}\subset\CQ_{\unl{d}}$ (parabolic quasiflags based at
$\infty\in\bC$) formed by the flags
$$0\subset\CW_1\subset\ldots\subset\CW_{n-1}\subset\CW=W\otimes\CO_\bC$$
such that $\CW_i\subset\CW$ is a vector subbundle in a
neighbourhood of $\infty\in\bC$, and the fiber of $\CW_i$ at
$\infty$ equals the span $\langle w_1,\ldots,w_{p_1+\ldots+p_i}\rangle\subset W$.

It is known to be a smooth connected quasiprojective variety of dimension
$\sum_{i=1}^{n-1}d_i(p_i+p_{i+1})$. Moreover, there is a natural proper morphism
$\fQ_{\unl{d}}\to \QM^{\unl{d}}_{G,P}$ and 
according to A.~Kuznetsov, this morphism is a small resolution of singularities
(cf.~Remark after~Theorem~7.3 of~\cite{BFGM}), so that
$H^\bullet_{T\times\BC^*}(\fQ_{\unl{d}})=
\on{IH}^\bullet_{T\times\BC^*}(\QM^{\unl{d}}_{G,P})$, and
$H^\bullet_{L\times\BC^*}(\fQ_{\unl{d}})=
\on{IH}^\bullet_{L\times\BC^*}(\QM^{\unl{d}}_{G,P})$.

\ssec{fixed points}{Fixed points}
The group $G\times\BC^*$ acts naturally on
$\CQ_{\unl{d}}$, and the group $T\times\BC^*$ acts
naturally on $\fQ_{\unl{d}}$. The set of fixed points of
$T\times\BC^*$ on $\fQ_{\unl{d}}$ is finite; its description is absolutely
similar to~\cite{FNR},~2.2, which we presently recall.

Let $\vec{\unl{d}}$ be a collection of nonnegative integral vectors
$\vec{d}_{ij}=(d_{ij}^{(1)},\ldots,d_{ij}^{(p_j)}),\ n-1\geq i\geq j\geq 1$,
such that $d_i=\sum_{j=1}^i|d_{ij}|=\sum_{j=1}^i\sum_{l=1}^{p_j}d_{ij}^{(l)}$,
and for $i\geq k\geq j$ we have $\vec{d}_{kj}\geq \vec{d}_{ij}$, i.e. for
any $1\leq l\leq p_j$ we have $d_{kj}^{(l)}\geq d_{ij}^{(l)}$.
Abusing notation we denote by $\vec{\unl{d}}$ the corresponding
$T\times\BC^*$-fixed point in $\fQ_{\unl{d}}$:

$$
\CW_1=\CO_\bC(-d_{11}^{(1)}\cdot0)w_1\oplus\ldots\oplus
\CO_\bC(-d_{11}^{(p_1)}\cdot0)w_{p_1},
$$

$$
\begin{aligned}
\CW_2=\CO_\bC(-d_{21}^{(1)}\cdot0)w_1\oplus\ldots\oplus
\CO_\bC(-d_{21}^{(p_1)}\cdot0)w_{p_1}\oplus\CO_\bC(-d_{22}^{(1)}\cdot0)w_{p_1+1}
\oplus\ldots\\
\ldots\oplus\CO_\bC(-d_{22}^{(p_2)}\cdot0)w_{p_1+p_2},
\end{aligned}
$$

$$\ldots\ \ldots\ \ldots\ ,$$

$$
\begin{aligned}
\CW_{n-1}=\CO_\bC(-d_{n-1,1}^{(1)}\cdot0)w_1\oplus\ldots\oplus
\CO_\bC(-d_{n-1,1}^{(p_1)}\cdot0)w_{p_1}\oplus\ldots\\
\ldots\oplus
\CO_\bC(-d_{n-1,n-1}^{(1)}\cdot0)w_{p_1+\ldots+p_{n-2}+1}
\oplus\ldots\oplus\CO_\bC(-d_{n-1,n-1}^{(p_{n-1})}\cdot0)w_{p_1+\ldots+p_{n-1}}.
\end{aligned}
$$

\ssec{classic}{Correspondences}
For $i\in\{1,\ldots,n-1\}$, and $\unl{d}=(d_1,\ldots,d_{n-1})$, we
set $\unl{d}+i:=(d_1,\ldots,d_i+1,\ldots,d_{n-1})$. We have a
correspondence $\CE_{\unl{d},i}\subset\CQ_{\unl{d}}\times
\CQ_{\unl{d}+i}$ formed by the pairs $(\CW_\bullet,\CW'_\bullet)$
such that for $j\ne i$ we have $\CW_j=\CW'_j$, and
$\CW'_i\subset\CW_i$, cf.~\cite{FNR},~2.3. In other words,
$\CE_{\unl{d},i}$ is the moduli space of flags of locally free
sheaves
$$0\subset\CW_1\subset\ldots\CW_{i-1}\subset\CW'_i\subset\CW_i\subset
\CW_{i+1}\ldots\subset\CW_{n-1}\subset\CW$$ such that
$\on{rank}(\CW_k)=p_1+\ldots+p_k$, and $\deg(\CW_k)=-d_k$, while
$\on{rank}(\CW'_i)=p_1+\ldots+p_i$, and $\deg(\CW'_i)=-d_i-1$.

According to~\cite{L1},~2.10, $\CE_{\unl{d},i}$ is a smooth
projective algebraic variety of dimension
$\dim(G/P)+\sum_{i=1}^{n-1}d_i(p_i+p_{i+1})+p_i$

We denote by $\bp$ (resp. $\bq$) the natural projection
$\CE_{\unl{d},i}\to\CQ_{\unl{d}}$ (resp.
$\CE_{\unl{d},i}\to\CQ_{\unl{d}+i}$). We also have a map $\br:\
\CE_{\unl{d},i}\to\bC,$
$$(0\subset\CW_1\subset\ldots\CW_{i-1}\subset\CW'_i\subset\CW_i\subset
\CW_{i+1}\ldots\subset\CW_{n-1}\subset\CW)\mapsto\on{supp}(\CW_i/\CW'_i).$$

The correspondence $\CE_{\unl{d},i}$ comes equipped with a natural
line bundle $\CL_i$ whose fiber at a point
$$(0\subset\CW_1\subset\ldots\CW_{i-1}\subset\CW'_i\subset\CW_i\subset
\CW_{i+1}\ldots\subset\CW_{n-1}\subset\CW)$$ equals
$\Gamma(\bC,\CW_i/\CW'_i)$. Let $q$ stand for the character of $T\times\BC^*:\
(\unl{t},v)\mapsto v$. We define the line bundle $\CL'_i:=q^{1-i}\CL_i$
on the correspondence $\CE_{\unl{d},i}$, that is $\CL'_i$ and $\CL_i$
are isomorphic as line bundles but the equivariant structure of $\CL'_i$
is obtained from the equivariant structure of $\CL_i$ by the twist by the
character $q^{1-i}$.

Finally, we have a transposed correspondence
$^\sT\CE_{\unl{d},i}\subset \CQ_{\unl{d}+i}\times\CQ_{\unl{d}}$.

\ssec{}{}
Restricting to $\fQ_{\unl{d}}\subset\CQ_{\unl{d}}$ we obtain the
correspondence
$\fE_{\unl{d},i}\subset\fQ_{\unl{d}}\times\fQ_{\unl{d}+i}$ together
with line bundle $\fL_i$ and the natural maps $\bp:\
\fE_{\unl{d},i}\to\fQ_{\unl{d}},\ \bq:\
\fE_{\unl{d},i}\to\fQ_{\unl{d}+i},\ \br:\
\fE_{\unl{d},i}\to\bC-\infty$. We also have a transposed
correspondence $^\sT\fE_{\unl{d},i}\subset
\fQ_{\unl{d}+i}\times\fQ_{\unl{d}}$. It is a smooth quasiprojective
variety of dimension $\sum_{i=1}^{n-1}d_i(p_i+p_{i+1})+p_i$.

\ssec{}{}
We denote by ${}'\IH_{G,P,T}$ the direct sum of equivariant (complexified)
cohomology:
${}'\IH_{G,P,T}=\oplus_{\unl{d}}H^\bullet_{T\times\BC^*}(\fQ_{\unl{d}})$.
It is a module over
$H^\bullet_{T\times\BC^*}(pt)=\BC[\ft\oplus\BC]=
\BC[x_1,\ldots,x_N,\hbar]$. Here $\ft\oplus\BC$ is the
Lie algebra of $T\times\BC^*$. We define $\hbar$ as
the positive generator of $H^2_{\BC^*}(pt,\BZ)$. Similarly, we define
$x_i\in H^2_{T}(pt,\BZ)$ in terms of the corresponding
one-parametric subgroup.
We define $\IH_{G,P,T}=\ {}'\IH_{G,P,T}\otimes_{H^\bullet_{T\times\BC^*}(pt)}
\on{Frac}(H^\bullet_{T\times\BC^*}(pt))$.

We have an evident grading
$$
\IH_{G,P,T}=\bigoplus\limits_{\unl{d}}\IH_{G,P,T}^{\unl{d}},\quad\text{where}\quad
\IH_{G,P,T}^{\unl{d}}=H^\bullet_{T\times\BC^*}(\fQ_{\unl{d}})
\otimes_{H^\bullet_{T\times\BC^*}(pt)}
\on{Frac}(H^\bullet_{T\times\BC^*}(pt)).
$$

According to the Thomason localization theorem,
restriction to the $T\times\BC^*$-fixed
point set induces an isomorphism
$$H^\bullet_{T\times\BC^*}(\fQ_{\unl{d}})
\otimes_{H^\bullet_{T\times\BC^*}(pt)}
\on{Frac}(H^\bullet_{T\times\BC^*}(pt))\to
H^\bullet_{T\times\BC^*}(\fQ_{\unl{d}}^{T\times\BC^*})
\otimes_{H^\bullet_{T\times\BC^*}(pt)}
\on{Frac}(H^\bullet_{T\times\BC^*}(pt))$$

The fundamental cycles $[\vec{\unl{d}}]$ of the
$T\times\BC^*$-fixed points $\vec{\unl{d}}$ (see~\refss{fixed points})
form a basis in $\oplus_{\unl{d}}H^\bullet_{T\times\BC^*}
(\fQ_{\unl{d}}^{T\times\BC^*})
\otimes_{H^\bullet_{T\times\BC^*}(pt)}\on{Frac}
(H^\bullet_{T\times\BC^*}(pt))$. The embedding of a point
$\vec{\unl{d}}$ into $\fQ_{\unl{d}}$ is a proper morphism, so the
direct image in the equivariant cohomology is well defined, and we will
denote by $[\vec{\unl{d}}]\in \IH_{G,P,T}^{\unl{d}}$ the direct image of the
fundamental cycle of the point $\vec{\unl{d}}$. The set
$\{[\vec{\unl{d}}]\}$ forms a basis of $\IH_{G,P,T}$.

\ssec{yang cor}{}
For any $0\leq i\leq n$ we will denote by $\unl{\CW}{}_i$ the tautological
$(p_1+\ldots+p_i)$-dimensional vector bundle on $\fQ_{\unl{d}}\times\bC$.
By the K\"unneth formula we have
$H^\bullet_{T\times\BC^*}(\fQ_{\unl{d}}\times\bC)=
H^\bullet_{T\times\BC^*}(\fQ_{\unl{d}})\otimes1\oplus
H^\bullet_{T\times\BC^*}(\fQ_{\unl{d}})\otimes\tau$ where
$\tau\in H^2_{\BC^*}(\bC)$ is the first Chern class of $\CO(1)$.
Under this decomposition, for the Chern class $c_j(\unl{\CW}{}_i)$ we have
$c_j(\unl{\CW}{}_i)=:c_j^{(j)}(\unl{\CW}{}_i)\otimes1+
c_j^{(j-1)}(\unl{\CW}{}_i)\otimes\tau$
where $c_j^{(j)}(\unl{\CW}{}_i)\in
H^{2j}_{T\times\BC^*}(\fQ_{\unl{d}})$, and
$c_j^{(j-1)}(\unl{\CW}{}_i)\in
H^{2j-2}_{T\times\BC^*}(\fQ_{\unl{d}})$.

For $0\leq m\leq n$ we introduce the generating series $\sA_m(u)$
with coefficients in the equivariant cohomology ring of
$\fQ_{\unl{d}}$ as follows:
$$\sA_m(u):=u^{p_1+\ldots+p_m}+\sum_{r=1}^{p_1+\ldots+p_m}(-\hbar)^{-r}
\left(c_r^{(r)}(\unl{\CW}{}_m)-\hbar
c_r^{(r-1)}(\unl{\CW}{}_m)\right) u^{p_1+\ldots+p_m-r}.
$$
In
particular, $\sA_0(u):=1$.

We also define the operators
\eq{ddvas}
\se_k^{(r+1+p_{k+1}-p_k)}:=
\bp_*(c_1(\CL'_k)^r\cdot\bq^*):\ \IH_{G,P,T}^{\unl{d}}\to
\IH_{G,P,T}^{\unl{d}-k},\ r\geq0
\end{equation}

\eq{ttris}
\ssf_k^{(r+1)}:=-\bq_*(c_1(\CL'_k)^r\cdot\bp^*):\ \IH_{G,P,T}^{\unl{d}}\to
\IH_{G,P,T}^{\unl{d}+k},\
r\geq0
\end{equation}

We consider the following generating series of operators on $\IH_{G,P,T}$:
\eq{raz}
\sd_k(u)=1+\sum_{s=1}^\infty\sd_k^{(s)}\hbar^{-s+1}u^{-s}:=
\sa_k(u+k-1)\sa_{k-1}(u+k-1)^{-1}:\ \IH_{G,P,T}^{\unl{d}}\to
\IH_{G,P,T}^{\unl{d}}[[u^{-1}]],
\end{equation}
where $ 1\leq k\leq n$ and
\eq{razz}
\sa_k(u):=u^{-p_1}(u-1)^{-p_2}\ldots(u-k+1)^{-p_k}\sA_k(u);
\end{equation}

\eq{dvas}
\se_k(u)=\sum_{s=1+p_{k+1}-p_k}^\infty\se_k^{(s)}\hbar^{-s+1}u^{-s}:\
\IH_{G,P,T}^{\unl{d}}\to \IH_{G,P,T}^{\unl{d}-k}[[u^{-1}]], 1\leq k\leq n-1;
\end{equation}

\eq{tris}
\ssf_k(u)=\sum_{s=1}^\infty\ssf_k^{(s)}\hbar^{-s+1}u^{-s}:\
\IH_{G,P,T}^{\unl{d}}\to \IH_{G,P,T}^{\unl{d}+k}[[u^{-1}]], 1\leq k\leq n-1.
\end{equation}

We also introduce the auxiliary series $\sB_k(u),\sC_k(u)$ by the
formulas~\refe{dvass},~\refe{triss}.

The following Theorem is a straightforward generalization of~Theorem~2.9
and the proof of~Theorem~2.12 of~\cite{FNR}, which are in turn its particular
case for $p_1=\ldots=p_n=1$.
\th{matrix coeff}
The matrix coefficients of the operators $\se_i^{(s)},\ssf_i^{(s)}$
in the fixed point base $\{[\vec{\unl{d}}]\}$ of $\IH_{G,P,T}$ are as follows:
$$\se_{i[\vec{\unl{d}},\vec{\unl{d}}{}']}^{(s)}=
\hbar^{-1}(\sfp_{ij}^{(a)}-i\hbar)^{s-1-p_{i+1}+p_i}
\prod_{\substack{k\leq i,\ b\leq p_k\\ (k,b)\ne(j,a)}}
(\sfp_{ij}^{(a)}-\sfp_{ik}^{(b)})^{-1}
\prod_{k\leq i+1}\prod_{b\leq p_k}(\sfp_{ij}^{(a)}-\sfp_{i+1,k}^{(b)}),$$
if $d_{ij}^{(a)}{}'=d_{ij}^{(a)}-1$ for certain $j\leq i$;
$$\ssf_{i[\vec{\unl{d}},\vec{\unl{d}}{}']}^{(s)}=
-\hbar^{-1}(\sfp_{ij}^{(a)}+(1-i)\hbar)^{s-1}
\prod_{\substack{k\leq i,\ b\leq p_k\\ (k,b)\ne(j,a)}}
(\sfp_{ij}^{(a)}-\sfp_{ik}^{(b)})^{-1}
\prod_{k\leq i-1}\prod_{b\leq p_k}(\sfp_{ij}^{(a)}-\sfp_{i-1,k}^{(b)}),$$
if $d_{ij}^{(a)}{}'=d_{ij}^{(a)}+1$ for certain $j\leq i$.
All the other matrix coefficients of $\se_i^{(s)},\ssf_i^{(s)}$ vanish.
Furthermore, the eigenvalue of $\sA_i(u)$ on $[\vec{\unl{d}}]$ equals
$$\prod_{j\leq i}\prod_{a\leq p_j}(u-\hbar^{-1}\sfp_{ij}^{(a)}).$$
\eth

\prop{isom} The isomorphism $\Psi:\ \IH_{G,P,T}\iso V,\ [\vec{\unl{d}}]\mapsto
(-1)^{|\unl{d}|}\hbar^{\sum_{i=1}^{n-1}d_ip_i}\xi_{\vec{\unl{d}}}$
intertwines the same named operators $\sd_i,\se_i,\ssf_i$, etc. In
particular, the operators $\sd_i,\se_i,\ssf_i$ defined
in~\refe{raz},~\refe{ddvas},~\refe{ttris}, turn $\IH_{G,P,T}$ into the
Gelfand-Tsetlin module over $U^\hbar(\fg,e)\otimes\CK'$. \eprop

\prf A straightforward comparison of~\reft{matrix coeff}
and~\refp{fmo}. \epr

\ssec{L vs T}{}
Now we return to the localized $L\times\BC^*$-equivariant cohomology
$\IH_{G,P}=\IH_{G,P,T}^{W_L}$. Note that the action of $W_L$ on $\IH_{G,P,T}$
is {\em semilinear} with respect to the structure of $\CK'$-module, and also
it commutes with the action of correspondences since both the correspondences
and the line bundles $\CL_i$ are equipped with the action of $L$.
Hence under the identification $\Psi:\ \IH_{G,P,T}\iso V$ the $W_L$-action
on $\IH_{G,P,T}$ goes to the $W_L$-action on $V$ introduced just
before~\refc{VM}. Combining~\refc{VM} with~\refp{isom} we arrive at the
following theorem proving~\refco{main}(1) in the case $\fg=\fgl_N$.

\th{ne}
The isomorphism $\Psi:\ \IH_{G,P,T}\iso V$ restricted to $W_L$-invariants
gives the same named isomorphism of $U^\hbar(\fgl_N,e)\otimes\CK$-modules
$\Psi:\ \IH_{G,P}\iso M(\fgl_N,e)$.
\eth

\sec{whit}{Whittaker vector and Shapovalov form}

\prop{whitta} If $1_{\unl{d}}$ stands for the unit cohomology class of
$\fQ_{\unl{d}}$, then $\Psi(1_{\unl{d}})=\fw_{\unl{d}}$.
\eprop

\prf
For $p_{i+1}-p_i+1\leq s<p_{i+1}$, we have $\se_i^{(s)}1_{\unl{d}}=0$ for
degree reasons (it would have had a negative degree). Similarly,
$\hbar\se_i^{(p_{i+1})}1_{\unl{d}+i}$, having degree 0, must be a constant
multiple of $1_{\unl{d}}$. More precisely, we decompose the projection
$\bp:\ \fE_{\unl{d},i}\to\fQ_{\unl{d}}$ into composition of the {\em proper}
$\bp\times\br:\ \fE_{\unl{d},i}\to\fQ_{\unl{d}}\times(\bC-\infty)$, and further
projection $\on{pr}:\ \fQ_{\unl{d}}\times(\bC-\infty)\to\fQ_{\unl{d}}$
with fibers $\bC-\infty=\BA^1$. We have $\bp_*(c_1(\CL_i)^{p_i-1}\cdot\bq^*
1_{\unl{d}+i})=\on{pr}_*(\bp\times\br)_*(c_1(\CL_i)^{p_i-1}\cdot\bq^*
1_{\unl{d}+i})$. Now $(\bp\times\br)_*(c_1(\CL_i)^{p_i-1}\cdot\bq^*1_{\unl{d}+i})$
is well defined in {\em nonlocalized} equivariant cohomology, and for the
degree reasons must take $1_{\unl{d}+i}$ to a constant multiple $c$ of the unit
class in the equivariant cohomology of $\fQ_{\unl{d}}\times(\bC-\infty)$.
Furthermore, $\on{pr}_*c=\hbar^{-1}c1_{\unl{d}}$. So it remains to calculate
the constant $c$. This can be done over the open subset $U\subset\fQ_{\unl{d}}$
where $\unl{\CW}{}_i/\unl{\CW}{}_{i-1}$ has no torsion, and hence
$\bp\times\br$ is a fibration with a fiber $\BP^{p_i-1}$. More precisely,
the correspondence $\fE_{\unl{d},i}$ over $U\times(\bC-\infty)$ is just the
projectivized vector bundle $\BP(\unl{\CW}{}_i/\unl{\CW}{}_{i-1})$, and
$\CL_i$ is nothing else than $\CO(1)$. We conclude that $c=1$.
The proposition is proved.
\epr

\prop{sha}
For $x,u\in\IH_{G,P}^{\unl{d}}$ we have
$(\Psi(x),\Psi(u))=(-1)^{|\unl{d}|}\int_{\fQ_{\unl{d}}}(xu)$.
\eprop

\prf
Evidently, the operators $\ssf_i^{(s)}$ and $-\se_i^{(s+p_{i+1}-p_i)}$ are
adjoint with respect to the pairing $\int(?\cdot?)$.
\epr

Thus we have fully proved \refco{main} for $\grg=\fgl(N)$.

%%%%%%%%%%%%%%%%%%%%%%%%%%%%%%%%%%%%%%%%%%%%%%%%%%%%%%%%%%%%%%%%%%%%%%%%%%%%%%%%%%%%%%%%%%%%%%%%%%%

\sec{AGT}{Relation to the AGT conjecture}
\ssec{}{The Uhlenbeck spaces of $\BA^2$}
Let $G$ be an almost simple simply connected complex algebraic group with maximal torus $T$ and let
$\grg,\grt$ be the corresponding Lie algebras.
For an integer $a\geq 0$ let $\Bun_G^d(\AA^2)$ denote the moduli space of principal $G$-bundles on
$\PP^2$ of second Chern class $a$ with a chosen trivialization at infinity (i.e. a trivialization
on the ``infinite" line $\PP^1_{\infty}$)
It is shown in \cite{bfg} that
 this space has the following properties:

 a) $\Bun_G^d(\AA^2)$ is non-empty if and only if $a\geq 0$;

 b) For $a\geq 0$ the space $\Bun_G^d(\AA^2)$ is an irreducible smooth
 quasi-affine variety of dimension $2a\check{h}$ where $\check{h}$ denotes the
dual Coxeter number of $G$.

In \cite{bfg} we construct an affine scheme $\calU^d_G(\AA^2)$ containing $\Bun_G^d(\AA^2)$
as a dense open subset which we are going to  call the {\it Uhlenbeck space} of bundles on $\AA^2$.

The scheme $\calU^d_G(\AA^2)$ is still irreducible but in general it is highly singular. The main property of $\calU^d_G(\AA^2)$
is that it possesses the following stratification:
\eq{stratification}
\Bun_G^d(\AA^2)=\bigcup\limits_{0\leq b\leq a}\Bun_G^b(\AA^2)\x \Sym^{a-b}(\AA^2).
\end{equation}
Here each $\Bun_G^b(\AA^2)\x \Sym^{a-b}(\AA^2)$ is a locally closed subset of $\calU^d_G(\AA^2)$ and its closure is
equal to the union of similar subsets corresponding to all $b'\leq b$.

We shall denote by $\Bun_G(\AA^2)$ (resp. $\calU_G(\AA^2)$) the disjoint union of all $\Bun_G^d(\AA^2)$
(resp. of $\calU^d_G(\AA^2)$).

Let us note that the group $G\x \GL(2)$ acts naturally on $\Bun_G^d(\AA^2)$: here the first factor acts by changing
the trivialization at ${\infty}$ and the second factor acts on $\AA^2$.
It is easy to deduce from the construction of \cite{bfg} that this action extends to an action of the same group
on the Uhlenbeck space $\calU^d_G(\AA^2)$.

The group $G\x \GL(2)$ acts naturally on $\calU^d_G$ where $G$ acts by
changing the trivialization
at $\PP^1_{\infty}$ and $\GL(2)$ acts on $\PP^2$ preserving $\PP^1_{\infty}$.

%--------------------------------------------------------------------------------------
\ssec{}{Instanton counting}We may now consider  the  equivariant integral
$$
\int\limits_{\calU^d_G}1^d
$$
of the unit $G\times \GL(2)$-equivariant cohomology class (which we
denote by $1^d$) over
$\calU_G^d$; the integral takes values in the field $\calK$ which is
the field of fractions of the algebra $\calA=H^*_{G\x \GL(2)}(pt)$.
Note that $\calA$ is
canonically isomorphic to the algebra of polynomial functions on the Lie algebra
$\grg\x\fgl(2)$  which
are invariant with respect to the adjoint action.
Thus each $\int\limits_{\calU^d_G}1^d$ may naturally be regarded as a
rational function of $a\in\grt$ and $(\eps_1,\eps_2)\in \CC^2$; this function
must be invariant with respect to the natural action of $W$ on $\grt$ and with respect
to interchanging $\eps_1$ and $\eps_2$.

Consider now the generating function
$$
\calZ=\sum\limits_{d=0}^\infty Q^d \int\limits_{\calU_G^d} 1^d.
$$
It can (and should) be thought of as a function of the variables
$\grq$ and $a,\eps_1,\eps_2$ as before. The function $Z(Q,a,\eps_1,\eps_2)$ is called the {\em Nekrasov partition
function} of pure $N=2$ supersymmetric gauge theory.

%----------------------------------------------------------------------------------------------------------------------
\ssec{agt-sl2}{The AGT conjecture}In \cite{AGT} Alday, Gaiotto and Tachikawa suggested
a relation between 4-dimensional supersymmetric gauge theory for $G=\SL(2)$ and the so called Liouville
2-dimensional conformal field theory; some generalizations to other groups were suggested in~\cite{AY},~\cite{MM}.
Here we are going to formulate a few mathematical statements suggested by the AGT conjecture (it is not clear
to us whether from the physics point of view they should be perceived
as direct corollaries of it).

Consider the above Uhlenbeck space $\calU^d_G$ and let
$\IH_{G\x \
GL(2)}(\calU^d_G)$ denote its equivariant intersection cohomology. This is a module over the algebra
$\calA_{G\x \GL(2)}:=H^*_{G\x \GL(2)}(pt)$; this algebra is just the algebra of polynomial functions
on $\grg\x \mathfrak{gl}(2)$ which are invariant under the adjoint action. We denote by $\calK_{G\x \GL(2)}$ its field
of fractions and we let
$$
\IH^{d,\aff}_G=\IH_{G\x \GL(2)}^*(\calU^d_G)\underset{\calA_{G\x \GL(2)}}
\ten\calK_{G\x \GL(2)}.
$$
We also set $\IH_G^{\aff}$ to be the direct sum of all the $\IH^{d,\aff}_G$.  This is a vector space
over $\calK_{G\x \GL(2)}$ which informally we may think of as a family of vector
spaces parametrized by $a\in\grt/W$ and $(\eps_1,\eps_2)\in \CC^2/\ZZ_2$.
Also, each $\IH^d_G$ is endowed with a perfect symmetric pairing $\la\cdot,\cdot\ra$, which is equal to
the Poincar\'e pairing multiplied by $(-1)^{\check{h}d}$.

Consider now the case $G=\SL(2)$. Then we can identify the Cartan subalgebra $\grh$ of $\mathfrak{sl}(2)$ with $\CC$.
Thus $a\in\grh$ can be thought of as a complex number.

\medskip
\noindent
{\bf Warning.} It is important to note that if we think about $a$ as a weight of the
Langlands dual $\mathfrak{sl}(2)$-algebra with standard generators $e,h,f$, then by definition the value of this
weight on $h$ is equal to $\alp(a)=2a$ (where $\alp$ denotes the simple root of $\mathfrak{sl}(2)$. This observation
will be used below.

\medskip
\noindent
Let $\Vir$ denote the Virasoro Lie algebra; it has the standard generators $\{ L_n\}_{n\in \ZZ}$ and $\bfc$
where $\bfc$ is central and $L_n$'s satisfy the standard relations.

Given a field $K$ of characteristic $0$, for every $\Del\in K,c\in K$ we may
consider the Verma module
$\bM_{\Del,c}$ over $\Vir$ on which $\bfc$ acts by $c$ and which is generated by
a vector $m_{\Del,c}$ such that
$$
L_0(m_{\Del,c})=\Del m_{\Del,c};\qquad L_n(m_{\Del,c})=0\quad\text{for $n>0$}.
$$

In addition, the Verma module ${\bf M}_{\Del,c}$ is equal to the direct sum of its $L_0$-eigen-spaces
${\bf M}_{\Del,c,d}$ where $n\in \ZZ_{\geq 0}$ and $L_0$ acts on ${\bf M}_{\Del,c,d}$ by $\Del+d$. It is easy to see that
there exists unique collection of vectors $w_{\Del,c,d}\in {\bf M}_{\Del,c,d}$ (for all $d\in\ZZ_{\geq 0}$) such that

1) $w_{\Del,c,0}=m_{\Del,c}$;

2) We have $L_i\cdot w_{\Del,c,d}=0$ for $i>1$ and $L_1\cdot w_{\Del,c,d}=w_{\Del,c,d-1}$.

We let $w_{\Del,c}$ denote the sum of all the $w_{\Del,c,d}$; this is an element of the completed
Verma module $\hat{\bf M}_{\Del,c}=\prod_{d\geq 0}{\bf M}_{\Del,c,d}$.

In addition the module ${\bf M}_{\Del,c}$ possesses unique symmetric bilinear
form $\la\cdot,\cdot\ra$ such
that $\la m_{\Del,c},m_{\Del,c}\ra=1$ and $L_n$ is adjoint to $L_{-n}$.

Then the AGT conjecture implies the following (the statement below is often
referred to as ``non-conformal limit"
of the AGT conjecture, cf. \cite{MMM}):
\conj{agt-sl2}
Let
\eq{agt-param}
\Del=-\frac{a^2}{\eps_1\eps_2}+\frac{(\eps_1+\eps_2)^2}{4\eps_1\eps_2};\quad c=1+\frac{6(\eps_1+\eps_2)^2}{\eps_1\eps_2}.
\end{equation}
Then
$$
\int\limits_{\calU^d_G}1=\la w_{\Del,c,d},w_{\Del,c,d}\ra.
$$
In other words,
$$
Z(a,\eps_1,\eps_2, (-1)^{\check{h}}Q)=Q^{-\Del}\la w_{\Del,c},Q^{L_0}w_{\Del,c}\ra.
$$
\econj
In fact, it is quite natural to expect that the following stronger result holds:
\conj{agt-sl2-conceptual}
\begin{enumerate}
\item
There exists an action of the Virasoro algebra $\Vir$ on $\IH_G$ such that
with this action $\IH_G^{\aff}$ becomes isomorphic to ${\bf M}_{\Del,c}$ where

\item The intersection pairing is $\Vir$-invariant, i.e. the adjoint operator
to $L_n$ is $L_{-n}$.
\item
$L_n\cdot 1^d=0$ for any $n>1$ and $d\geq 0$.
\item
$L_1 \cdot 1^d=1^{d-1}$ for any $d>0$.
\end{enumerate}
\econj

One can generalize \refco{agt-sl2} and \refco{agt-sl2-conceptual}
to arbitrary $G$. We are not going to give details here,
but let us stress one thing: when $G$ is simply laced the Virasoro
algebra $\Vir$ has to be replaced
by the {\em $W$-algebra} corresponding to the affine Lie algebra $\grg_{\aff}$; in fact for general
$G$ (not necessarily simply laced) we believe that the $W$-algebra associated with the Langlands
dual affine Lie algebra $\grg_{\aff}^{\vee}$ should appear (cf. the next subsection for
some motivation). In other words, we expect that for general $G$ the space
$\IH_G$ carries a natural action of the $W$-algebra $W(\grg_{\aff}^{\vee})$, which makes it into a Verma
module over this algebra; analogs of properties~2,3,4
above are also expected\footnote{A modification of this conjecture exists also for
$\grg=\fgl(N)$; this conjecture will
be proved in \cite{MO}.}.
It is easy to deduce from the results of \cite{bfg}  that the character of $\IH_G$ is equal to the
character of the Verma module for the $W$-algebra.

Before we explain the connection of \refco{agt-sl2-conceptual} with the rest of the paper, let us recall some
(known) modification of it.
%---------------------------------------------------------------------------------------------------------
\ssec{}{The flag case} Choose a parabolic
subgroup $P\subset G$ and let $\bfC$ denote the ``horizontal" line in $\PP^2$
(i.e. we choose a straight
line in $\PP^2$ different from the one at infinity and call).
Let $\Bun_{G,P}$ denote the moduli space of
the following objects:

1) A principal $G$-bundle $\calF_G$ on $\PP^2$;

2) A trivialization of $\calF_G$ on $\PP^1_{\infty}$;

3) A reduction of $\calF_G$ to $P$ on $\bfC$ compatible with the
trivialization of $\calF_G$ on $\bfC$.

Let us describe the connected components of $\Bun_{G,P}$. We are going to use the notations
of \refss{setup}. Let also
$\Lam_{G,P}^{\aff}=\Lam_{G,P}\x \ZZ$ be the lattice of characters
of $Z(\check{L})\x\CC^*$. Note that $\Lam_{G,G}^{\aff}=\ZZ$.

The lattice $\Lam^{\aff}_{G,P}$ contains canonical semi-group
$\Lam^{\aff,+}_{G,P}$ of positive elements (cf. \cite{bfg}). It is not difficult to see that the connected
components of $\Bun_{G,P}$ are parameterized  by the elements of
$\Lam_{G,P}^{\aff,+}$:
$$
    \Bun_{G,P}=\bigcup\limits_{\theta_{\aff}\in\Lam_{G,P}^{\aff,+}}
    \Bun_{G,P}^{\theta_{\aff}}.
$$

Typically, for $\theta_{\aff}\in \Lam_{G,P}^{\aff}$ we shall write
$\theta_{\aff}=(\theta,d)$ where $\theta\in \Lam_{G,P}$ and $d\in \ZZ$.

Each $\Bun_{G,P}^{\theta_{\aff}}$ is naturally acted upon by $P\x(\CC^*)^2$;
embedding $L$ into $P$ we get an action of $L\x(\CC^*)^2$ on
$\Bun_{G,P}^{\theta_{\aff}}$. In \cite{bfg} we define for each
$\theta_{\aff}\in\Lam_{G,P}^{\aff,+}$ certain Uhlenbeck scheme
$\calU_{G,P}^{\theta_{\aff}}$ which contains $\Bun_{G,P}^{\theta_{\aff}}$ as a dense
open subset. The scheme $\calU_{G,P}^{\theta_{\aff}}$ still admits an
action of $L\x(\CC^*)^2$.

Following \cite{B} define
\eq{partition-affine}
    \calZ_{G,P}^{\aff}=\sum\limits_{\theta\in\Lam_{G,P}^{\aff}} \grq_{\aff}^{\theta_{\aff}}
    \ \int\limits_{\calU_{G,P}^{\theta_{\aff}}}1_{G,P}^{\theta_{\aff}}.
\end{equation}

\medskip
\noindent
{\bf Remark.} In addition to \cite{B} and \cite{BrEt} various examples of functions $\calZ_{G,P}^{\aff}$ were studied
recently in the physical literature (cf. for example \cite{AGGTV}, \cite{AT}) as (the instanton part of)
the Nekrasov partition function in the presence of surface operators.

\medskip
One should think of $\calZ_{G,P}^{\aff}$ as a formal power series
in $\grq_{\aff}\in Z(\chL)\x\CC^*$ with values in the space of ad-invariant
rational functions on $\grl\x\CC^2$. Typically, we shall write
$\grq_{\aff}=(\grq,Q)$ where $\grq\in Z(\check{L})$ and $Q\in\CC^*$.
Also we shall denote an element of $\grl\x\CC^2$ by
$(a,\eps_1,\eps_2)$ or (sometimes it will be more convenient) by
$(a,\hbar,\eps)$ (note that for general $P$ (unlike in the case $P=G$)
the function $\calZ_{G,P}^{\aff}$ is not symmetric with respect to
switching
$\eps_1$ and $\eps_2$).

As before, let us now denote by $\IH_{G,P}^{\theta_{\aff}}$ the localized $L\x(\CC^*)^2$-equivariant
intersection cohomology of $\calU_{G,P}^{\theta_{\aff}}$; we also set $\IH_{G,P}^{\aff}$ to be the direct
sum of all the $\IH_{G,P}^{\theta_{\aff}}$; note that $\IH_{G,G}^{\aff}=\IH_{G}^{\aff}$. Then $\IH_{G,P}^{\aff}$ is $\Lam_{G,P}^{\aff,+}$-graded
$\calK_{L\x(\CC^*)^2}$ vector space. We can endow with a non-degenerate pairing $\la\cdot,\cdot\ra$

The following result is proved in \cite{B}:
\th{affine-flag}
Let $P=B$ be a Borel subgroup. Then $\IH_{G,B}^{\aff}$ possesses a natural action of the Lie
algebra $\check{\grg}_{\aff}$ such that
\begin{enumerate}
\item
As a $\check{\grg}_{\aff}$-module $\IH_{G,B}^{\aff}$ is isomorphic to $M(\lam_{\aff})$, where
$\lam_{\aff}=-\frac{(a,\eps_1)}{\eps_2}-\rho_{\aff}$.
\item
$\IH_{G,B}^{\theta_{\aff}}$ is the
$-\frac{(a,\eps_1)}{\eps_2}-\rho_{\aff}-\theta_{\aff}$-weight space of $\IH_{G,B}^{\aff}$.
\footnote{Here by $(a,\eps_1)$ we mean the weight of $\grg_{\aff}^{\vee}$
whose ``finite" component is $a$ and whose central
charge is $\eps_1$. Also $\rho_{\aff}$ is a weight of $\grg_{\aff}^{\vee}$ which takes value 1 on every simple coroot.}
\item The isomorphism of (1) takes the pairing $\langle\cdot,\cdot\rangle$
on $\IH_{G,B}^{\theta_{\aff}}$ to $(-1)^{|\theta_{\aff}|}$ times the Shapovalov
pairing $(\cdot,\cdot)$ on the
$-\frac{(a,\eps_1)}{\eps_2}-\rho_{\aff}-\theta_{\aff}$-weight space of
$M(\lam_{\aff})$.
\end{enumerate}
\eth

%-----------------------------------------------------------
\ssec{}{Interpretation via maps and the ``finite-dimensional" analog}
Choose now another smooth
projective curve $\bfX$ of genus $0$ and with two marked points
$0_\bfX,\infty_\bfX$. Choose also a coordinate $x$ on $\bfX$ such
that $x(0_\bfX)=0$ and $x(\infty_\bfX)=0$. Let us denote by
$\CG_{G,P,\bX}$ the scheme classifying triples
$(\F_G,\beta,\gamma)$, where

1) $\F_G$ is a principal $G$-bundle on $\bfX$;

2) $\beta$ is a trivialization of $\F_G$ on the formal
neighborhood of $\infty_\bfX$;

3)  $\gamma$ is  a reduction to $P$ of the fiber of $\F_G$ at
$0_\bX$.

We shall usually omit $\bfX$ from the
notations.
We shall also write $\calG_{G}^{\aff}$ for
$\calG_{G,G}^{\aff}$.

Let $e_{G,P}^{\aff}\in\calG_{G,P}^{\aff}$ denote the point
corresponding to the trivial $\calF_G$ with the natural $\beta$
and $\gamma$. It is explained in \cite{bfg} that the variety
$\Bun_{G,P}$ is canonically isomorphic to the scheme classifying
{\it based maps} from $(\PP^1,\infty)$ to
$(\calG_{G,P}^{\aff},e_{G,P}^{\aff})$ (i.e. maps from $\PP^1$
to $\calG_{G,P}$ sending $\infty$ to
$e_{G,P}^{\aff}$).

The scheme $\calG_{G,P}^{\aff}$ may (and should) be thought of as
a partial flag variety for $\grg_{\aff}$. Thus the scheme $\Bun_{G,P}$ should
be thought of as an affine analog of $\calM_{G,P}$. Also the flag Uhlenbeck scheme
$\calU_{G,P}$ should be thought of as an affine analog of the
scheme $\QM_{G,P}$ (this analogy is explained in more detail in \cite{bfg}).

Thus \reft{affine-flag}
can be considered as an affine version of \reft{flag} and \refco{agt-sl2-conceptual}
(together with its generalization to arbitrary $\grg$ mentioned above) is an affine version
of \refco{main} (in fact, there should be a more general version of this conjecture, dealing
not only with $\IH_G^{\aff}$ but with arbitrary $\IH_{G,P}^{\aff}$).
To conclude the paper we are going to  explain how to use this analogy with \refco{main} in order to derive
the formulas \refe{agt-param} from \reft{affine-flag}.
For $\chi,k\in\CC$ let $M(\chi,k)$ denote the Verma module over $\mathfrak{sl}(2)_{\aff}$ with central charge
$k$ and highest weight $\chi$ (i.e. the standard generator $h$ of $\mathfrak{sl}(2)$ acts on the highest weight
vector as multiplication by $\chi$). According to \cite{FF} the algebra ${\bf Vir}$ is obtained by certain
BRST reeduction from $U(\mathfrak{sl}(2)_{\aff})$; the corresponding BRST reduction of $M(\chi,k)$ is
equal to $\bfM_{\Del,c}$ where
\eq{ff}
\Del=\frac{(\chi+1)^2-(k+1)^2}{4(k+2)};\qquad c=1-\frac{6(k+1)^2}{k+2}.
\end{equation}
According to \reft{affine-flag} we should take
\eq{eqchic}
\chi=-\frac{2a}{\eps_2}-1;\footnote{The appearance of the factor 2 is
explained in the Warning in \refss{agt-sl2}.}
\qquad k=-\frac{\eps_1}{\eps_2}-2.
\end{equation}
Thus $k+2=-\eps_1/\eps_2$ and
$$
\frac{(k+1)^2}{k+2}=-\frac{(\eps_1+\eps_2)^2}{\eps_1\eps_2}.
$$
Hence
$$
c=1+\frac{6(\eps_1+\eps_2)^2}{\eps_1\eps_2}
$$
and
$$
\Del=\frac{(\chi+1)^2-(k+1)^2}{4(k+2)}=-\frac{a^2}{\eps_1\eps_2}+
\frac{(\eps_1+\eps_2)^2}{4\eps_1\eps_2},
$$
which coincides with~\refe{agt-param}.

%------------------------------------------------------

\bigskip
\footnotesize{
{\bf A.B.}: Department of Mathematics, Brown University,
151 Thayer St., Providence RI
02912, USA;\\
{\tt braval@math.brown.edu}}

\footnotesize{
{\bf B.F.}: Landau Institute and State University Higher School of Economics\\
Department of Mathematics, 20 Myasnitskaya st, Moscow 101000 Russia;\\
{\tt bfeigin@gmail.com}}

\footnotesize{
{\bf M.F.}: IMU, IITP and State University Higher School of Economics\\
Department of Mathematics, 20 Myasnitskaya st, Moscow 101000 Russia;\\
{\tt fnklberg@gmail.com}}

\footnotesize{
{\bf L.R.}: IITP and State University Higher School of Economics\\
Department of Mathematics, 20 Myasnitskaya st, Moscow 101000 Russia;\\
{\tt leo.rybnikov@gmail.com}}


\begin{thebibliography}{jjjj}

\bibitem{AGGTV}
L.~F.~Alday, D.~Gaiotto, S.~Gukov, Y.~Tachikawa, H.~Verlinde,
{\em Loop and surface operators in N=2 gauge theory and Liouville modular geometry}, arXiv:0909.0945

\bibitem{AGT}
L.~F.~Alday, D.~Gaiotto and Y.~Tachikawa, {\em Liouville correlation
functions from four-dimensional gauge theories},  Lett. Math. Phys. {\bf 91}
(2010),  no. 2, 167--197.

\bibitem{AT}
L.~F.~Alday, Y.~Tachikawa, {\em Affine SL(2) conformal blocks from 4d gauge theories},
arXiv:1005.4469

\bibitem{AY} H.~Awata, Y.~Yamada, {\em Five-dimensional AGT Relation and
the deformed $\beta$-ensemble}, preprint arXiv:1004.5122.

\bibitem{B} A.~Braverman, {\em Instanton counting via affine Lie algebras.
I. Equivariant $J$-functions of (affine) flag manifolds and Whittaker vectors},
Algebraic structures and moduli spaces, CRM Proc. Lecture Notes {\bf 38}
Amer. Math. Soc., Providence, RI (2004), 113--132.

\bibitem{Bicm}A.~Braverman,
{\em Spaces of quasi-maps and their applications},
International Congress of Mathematicians. Vol. II,
Eur. Math. Soc., Z\"urich, (2006), 1145-1170.

\bibitem{BrEt}
A.~Braverman and P.~Etingof,
{Instanton counting via affine Lie algebras II:
from Whittaker vectors to the Seiberg-Witten prepotential},
Studies in Lie theory, 61--78, Progr. Math., {\bf 243}, Birkh\"auser Boston, Boston, MA, 2006.


\bibitem{bfg}
A.~Braverman, M.~Finkelberg and D.~Gaitsgory,
{\em Uhlenbeck spaces via affine Lie algebras},
in: The unity of mathematics
(volume dedicated to I.~M.~Gelfand's 90th birthday),
 Progr. Math., {\bf 244}, Birkh{\"a}user Boston (2006), 17-135.


\bibitem{BFGM} A.~Braverman, M.~Finkelberg, D.~Gaitsgory, I.~Mirkovi\'c,
{\em Intersection cohomology of Drinfeld's compactifications.}
Selecta Math. (N.S.)  {\bf 8}  (2002),  no. 3, 381--418.

\bibitem{BG} J.~Brundan, S.~M.~Goodwin, {\em Good grading polytopes},
Proc. London Math. Soc. {\bf 94} (2007), 155-180.

\bibitem{BGK} J.~Brundan, S.~Goodwin, A.~Kleshchev, {\em Highest weight
theory for finite $W$-algebras}, Int. Math. Res. Not. (2008) no. 15,
Art. ID rnn051, 53 pp.

\bibitem{BK} J.~Brundan, A.~Kleshchev, {\em Representations of shifted
Yangians and finite $W$-algebras}, Mem. Amer. Math. Soc. {\bf 196}, no. 918
(2008), viii+107pp.

\bibitem{etingof}
P.~Etingof,
{\em Whittaker functions on quantum groups and $q$-deformed Toda
operators}, in:
Differential topology, infinite-dimensional Lie algebras, and
applications, Amer. Math. Soc. Transl. Ser. 2, {\bf 194}, Amer. Math. Soc.,
Providence, RI,
(1999), 9--25.

\bibitem{FF}
B.~Feigin and E.~Frenkel, {\em Representations of affine Kac-Moody algebras,
bosonization and resolutions},
Lett. Math. Phys. {\bf 19} (1990), 307--317.

\bibitem{FM}
M.~Finkelberg and I.~Mirkovi\'c,
{\em Semi-infinite flags. I. Case of global curve $\PP^1$}, in:
Differential topology, infinite-dimensional Lie algebras, and
applications, Amer. Math. Soc. Transl. Ser. 2, {\bf 194}, Amer. Math. Soc.,
Providence, RI,
(1999), 81--112.


\bibitem{FFKM}
Feigin, B., Finkelberg, M., Kuznetsov, A., Mirkovi\'c, I.,
{\em Semi-infinite flags. II. Local and global intersection cohomology of quasimaps' spaces}, in:
Differential topology, infinite-dimensional Lie algebras, and
applications, Amer. Math. Soc. Transl. Ser. 2, {\bf 194}, Amer. Math. Soc.,
Providence, RI,
(1999), 113--148.


\bibitem{FNR} B.~Feigin, M.~Finkelberg, A.~Negut, L.~Rybnikov,
{\em Yangians and cohomology rings of Laumon spaces}, preprint arXiv:0812.4656.

\bibitem{FMO} V.~Futorny, A.~Molev, S.~Ovsienko, {\em Gelfand-Tsetlin bases
for representations of finite $W$-algebras and shifted Yangians},
in ``Lie theory and its applications in physics VII", (H.~D.~Doebner
and V.~K.~Dobrev, Eds), Proceedings of the VII International Workshop, Varna, Bulgaria, June 2007. Heron
Press, Sofia (2008), 352--363, see also arXiv:0711.0552.

\bibitem{GK}
A.~Givental, B.~Kim, {\em Quantum cohomology of flag manifolds and Toda
lattices},  Comm. Math. Phys. {\bf 168}  (1995),  no. 3, 609--641.


\bibitem{Kim} B.~Kim, {\em Quantum cohomology of flag manifolds $G/B$
and quantum Toda lattices}, Ann. of Math. (2) {\bf 149} (1999), 129--148.

\bibitem{L1} G.~Laumon, {\em Un Analogue Global du C\^one Nilpotent},
Duke Math. Journal {\bf 57} (1988), 647--671.

\bibitem{L2} G.~Laumon, {\em Faisceaux Automorphes Li\'es aux S\'eries
d'Eisenstein}, Perspect. Math. {\bf 10} (1990), 227--281.

\bibitem{MM}
 A. Mironov, A. Morozov, {\em On AGT relation in the case of U(3)},
Nucl.Phys.B {\bf 825} (2010), 1-37.

\bibitem{MMM}
A.~Marshakov, A.~Mironov and A.~Morozov,
{\em On non-conformal limit of the AGT relations},
Phys. Lett. B {\bf 682} (2009), no. 1, 125--129.

\bibitem{MO}
D.~Maulik and A.~Okounkov, in preparation.

\bibitem{Taki}
M.~Taki, {\em On AGT Conjecture for Pure Super Yang-Mills and W-algebra},
preprint arXiv:0912.4789.




\end{thebibliography}
\end{document}